\documentclass[12pt]{article}

\usepackage{epsfig}

\newenvironment{proof}{{\it Proof:\/}}{$\Box$\vskip 0.08in}
\usepackage{amssymb}

 \newtheorem{theorem}{Theorem}[section]
 \newtheorem{lemma}[theorem]{Lemma}

 \newtheorem{definition}[theorem]{Definition}

\newtheorem{conjecture}[theorem]{Conjecture}

\begin{document}
\pagestyle{myheadings}
 
\thispagestyle{empty}
\markboth {{\sc j.h.przytycki}}
{{\sc Links with the same Jones type polynomials}}
\begin{center}
\begin{LARGE}
\baselineskip=10pt
{\bf Search for different links with the same  Jones' type  polynomials:  
Ideas from graph theory and statistical mechanics.}
\end{LARGE}
\\
\ \\
  J\'OZEF H.~PRZYTYCKI
\end{center}
\begin{quotation}
ABSTRACT.
\baselineskip=10pt
We describe, in this talk,\footnote{This is the detailed version of 
the talk given  at the Banach Center Colloquium on 24th March 1994 
(``W poszukiwaniu
nietrywialnego w\c ez\l a z trywialnym wielomianem Jonesa: grafy
i mechanika statystyczna").} 
three methods of constructing different links with the same
Jones type invariant. All three can be thought as generalizations of mutation.
The first combines the satellite construction with mutation. The second
uses the notion of rotant, taken from the graph theory, the third, invented
by Jones, transplants into knot theory the idea of the Yang-Baxter equation 
with the spectral parameter (idea employed by Baxter  in the theory of solvable models in statistical mechanics). We extend the Jones result and
relate it  to Traczyk's work on rotors of links. We also show further 
applications of the Jones idea, e.g. to 3-string links in the solid torus.
We stress the fact that ideas coming from various areas of mathematics (and
theoretical physics) has been fruitfully used in knot theory, and vice versa.
\\
\end{quotation}
\ \\
{\Large \bf 0\ \  Introduction}\\ \ \\
Exactly ten years ago, at spring of 1984, Vaughan Jones introduced his
(Laurent) polynomial invariant of links, $V_L(t)$.
He checked immediately that it distinguishes many knots which were not
taken apart by the Alexander polynomial, e.g. the right handed trefoil
knot from the left handed trefoil knot, and the square knot from
the granny knot; Fig. 0.1.

\par\vspace{1cm}
\begin{center}
\begin{tabular}{c}
\includegraphics[trim=0mm 0mm 0mm 0mm, width=.65\linewidth]
{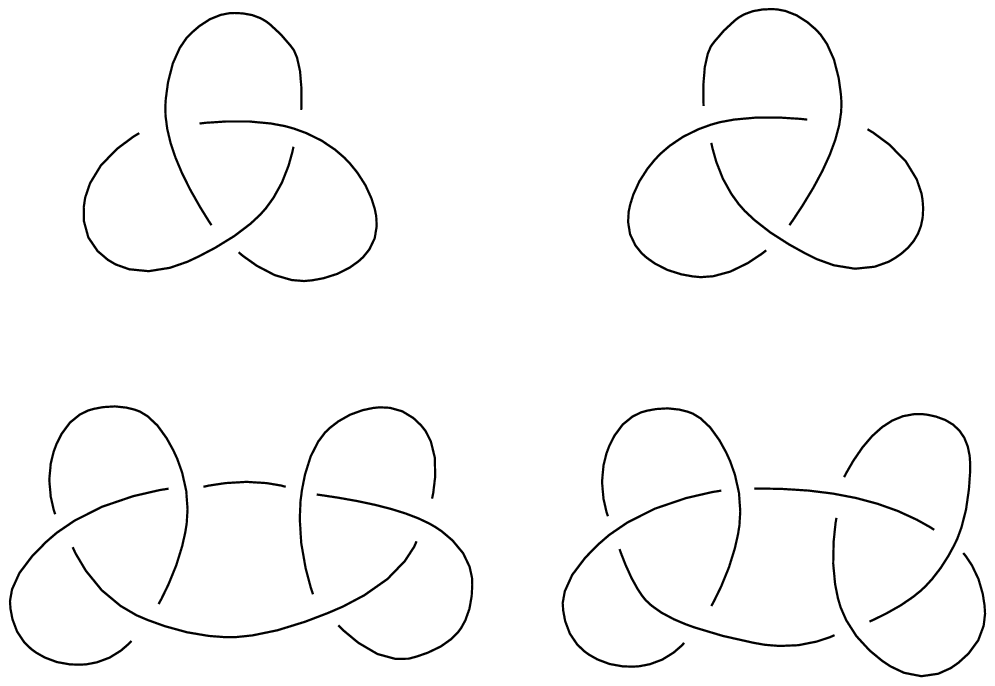}
\end{tabular}
\par\vspace{0.5cm}
Fig. 0.1
\end{center}

Jones also noticed that his polynomial is not universal. That is, there
are  different knots with the same polynomial; e.g. the Conway and
Kinoshita-Terasaka knots; Fig. 0.2. 
\par\vspace{1cm}
\begin{center}
\begin{tabular}{cc}
\includegraphics[trim=0mm 0mm 0mm 0mm, width=.3\linewidth]
{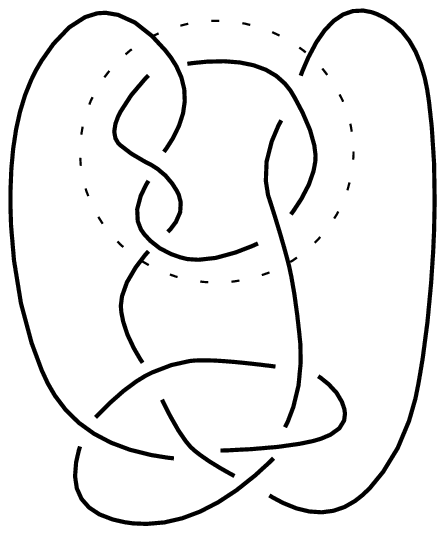}
&\includegraphics[trim=0mm 0mm 0mm 0mm, width=.3\linewidth]
{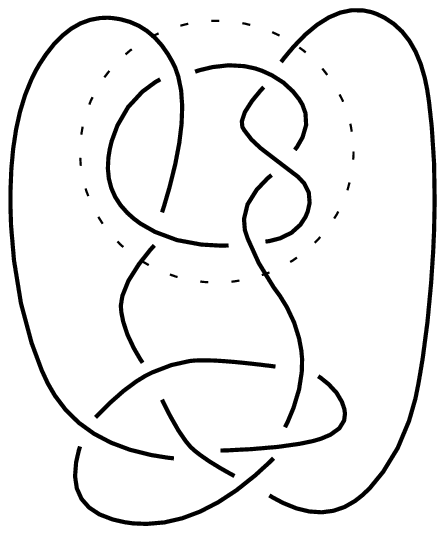}
\\
\end{tabular}
\par\vspace{0.5cm}
Fig. 0.2
\end{center}

Then Jones asked the fundamental
question whether there exists a nontrivial knot with the trivial 
polynomial. Ten year later this is still an open problem and specialists
differ in their opinion whether the answer is yes or no. In this talk,
I will concentrate on more accessible problem: how to construct different
links with the same Jones polynomial. It may shed some light into the
Jones question.\\ 
I will describe three methods of constructing knots with a coinciding
Jones polynomial (and its generalizations), each of which can be thought
as a generalization of the Conway idea of {\em mutation}.
First however, in the introduction, we remind the definitions of 
Jones type polynomials and of the Conway's mutation.
\begin{enumerate}
\item [(1)]
In the first part we consider satellites of mutants and their Jones
type invariants.
\item [(2)]
In the second part
we explain the idea of {\em rotors}.
\item [(3)]
In the third part we explore the idea of Jones of the spectral parameter
tangle.
\item [(4)]
In the fourth part we apply the idea to the skein polynomial of links.
\item [(5)]
In the fifth part we use a 3-string spectral parameter tangle.
\end{enumerate}

We remind now definitions of the Jones polynomial, $V_L(t)$, and its
generalizations: the skein (Homflypt) polynomial,  $P_L(v,z)$, 
and Kauffman polynomial, $F_L(a,x)$.
\begin{definition}\label{0.1}
The skein polynomial invariant of oriented links can be characterized by
the recursive relation (skein relation):
\begin{enumerate}
\item
[(i)] $v^{-1}P_{L_+}(v,z)-vP_{L_-}(v,z)=zP_{L_0}(v,z)$,\ where $L_+,L_-$ and
$L_0$ are three oriented link diagrams, which are the same outside a small
disk in which they look as in Fig. 0.3,
\end{enumerate}
and the initial condition
\begin{enumerate}
\item
[(ii)] $P_{T_1}=1$, \ where $T_1$ denotes the trivial
knot.
\end{enumerate}
\end{definition}

\par\vspace{1cm}
\begin{center}
\begin{tabular}{ccc}
\includegraphics[trim=0mm 0mm 0mm 0mm, width=.2\linewidth]
{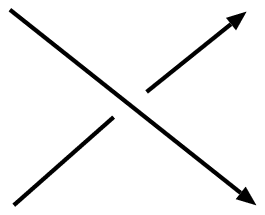}
&
\includegraphics[trim=0mm 0mm 0mm 0mm, width=.2\linewidth]
{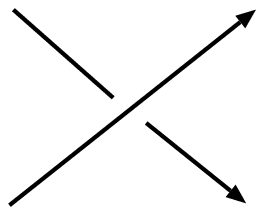}
&
\includegraphics[trim=0mm 0mm 0mm 0mm, width=.2\linewidth]
{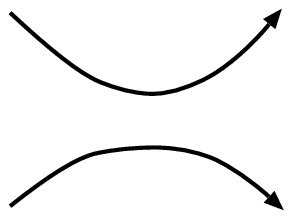}\\
$L_+$ & $L_-$ & $L_0$
\end{tabular}
\par\vspace{0.5cm}
Fig. 0.3
\end{center}

The Jones polynomial is defined as
$V_L(t)=P_L(t,\sqrt t-\frac{1}{\sqrt t})$.  The Alexander polynomial,  
${\Delta}_L(t)$, as normalized by Conway, satisfies \ ${\Delta}_L(t)=
P_L(1,\sqrt t-\frac{1}{\sqrt t})$.

Kauffman \cite{Ka-1} gave different approach to the Jones polynomial,
starting from the invariant of regular isotopy of unoriented
diagrams  or, equivalently, working with unoriented framed
links. This variant of the Jones polynomial is called now the Kauffman 
bracket polynomial.
\begin{definition}
The Kauffman bracket polynomial, $<L>\in Z[A^{\pm 1}]$, of framed 
unoriented links is characterized by the recursive relation (Kauffman
 bracket skein relation):
\begin{enumerate}
\item
[(i)] $<L_A>=A<L_0>+A^{-1}<L_{\infty}>$,\ where $L_A,L_0$ and $L_{\infty}$
denote diagrams of unoriented framed links, which are the same outside a
small disc in which they look as in Fig. 0.4. We use the convention that 
the framing of the diagram is vertical to the plane of projection, 
unless otherwise stated.
\item
[(ii)] Initial conditions:\\
$<T_n>=(-A^2-A^{-2})^{n-1}$,\ where $T_n$ is the trivial
framed n-component link.
\end{enumerate}
\end{definition}

\par\vspace{1cm}
\begin{center}
\begin{tabular}{cccc}
\includegraphics[trim=0mm 0mm 0mm 0mm, width=.17\linewidth]
{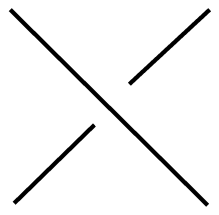}
&
\includegraphics[trim=0mm 0mm 0mm 0mm, width=.17\linewidth]
{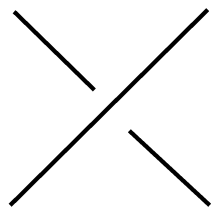}
&
\includegraphics[trim=0mm 0mm 0mm 0mm, width=.17\linewidth]
{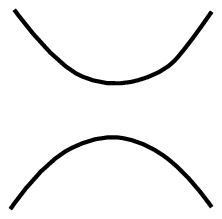}
&
\includegraphics[trim=0mm 0mm 0mm 0mm, width=.17\linewidth]
{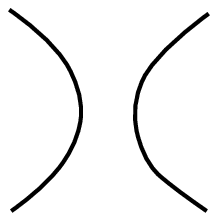}
\\
$L_A$ & $L_B$ & $L_0$ & $L_{\infty}$
\\
\end{tabular}
\par\vspace{0.5cm}
Fig. 0.4
\end{center}

If we orient $L$, then we get the Jones polynomial
$$V_L(A^{-4})=(-A^3)^{-Tait(L)}<L>,$$ \ where $Tait(L)$ is the sum of signs of
crossings of $L$. Equivalently, for a framed oriented link $L$, $Tait(L)$ is
``the defect" of the framing, that is the number of negative twists 
minus the number of positive twists which have 
to be performed on the framing of $L$ so that the new framing
agrees with that given by the Seifert surface of $L$. Let us use the
convention that $L^{(1)}$ denote a framed link obtained from $L$ by
twisting the framing of $L$ once in the positive (right handed) direction.
Notice that $L^{(1)}$ is not  uniquely defined if $L$ is not a knot, 
but its Kauffman bracket is well defined and the condition (ii) can be
replaced by:
\begin {enumerate}
\item
[(iii)] $<L^{(1)}> = -A^3<L>$.
\end{enumerate}

To introduce the Kauffman polynomial, it is, as before, very convenient to 
define it first for unoriented framed links.
\begin{definition}\label{0.3}
The Kauffman polynomial of framed unoriented links, ${\Lambda}_L(a,x)\in
 Z[a^{\pm 1},x^{\pm 1}]$,
 is characterized by the recursive relation (Kauffman
  skein relation):
\begin{enumerate}
\item
[(i)] $ {\Lambda}_{L_A}(a,x) + {\Lambda}_{L_B}(a,x) =
x({\Lambda}_{L_0}(a,x) + {\Lambda}_{L_{\infty}}(a,x))$, \
where $L_A,L_B,$ \ \ $L_0$ and $L_{\infty}$ are four
 diagrams of unoriented framed links, which are the same outside a
small disc in which they look as in Fig. 0.4.
\item
[(ii)] Framing relations:\\
$L^{(1)}=aL$ .
\item
[(iii)] Initial condition:\\
${\Lambda}_{T_1}=1$, for the trivial framed knot $T_1$.
\end{enumerate}

If we orient $L$, then the Kauffman polynomial of the oriented unframed
link $L$, $F_L(a,x)$, is defined to be:\\
$F_L(a,x)=a^{-Tait(L)}{\Lambda}_L(a,x)$.
\end{definition}

The simplest method of producing different links
with the same invariant is {\it mutation} invented before 1960 by J.Conway.
We describe the Conway idea of tangles and mutation below:\footnote{
Conway recalls working out a good part of his theory of tangles while still
a high school student \cite{Alb}.} \\
Consider a 2-tangle, $L$, that is a part of a link diagram placed in a disk,
with 4 boundary points (2 inputs and 2 outputs, if L is oriented); 
see Fig. 0.5.

\par\vspace{1cm}
\begin{center}
\begin{tabular}{ccc}
\includegraphics[trim=0mm 0mm 0mm 0mm, width=.25\linewidth]
{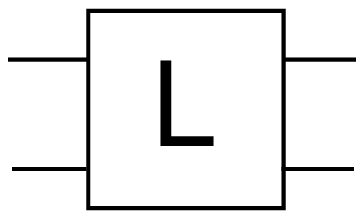}
&
\includegraphics[trim=0mm 0mm 0mm 0mm, width=.25\linewidth]
{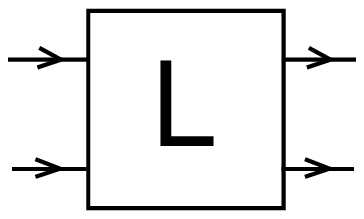}
&
\includegraphics[trim=0mm 0mm 0mm 0mm, width=.25\linewidth]
{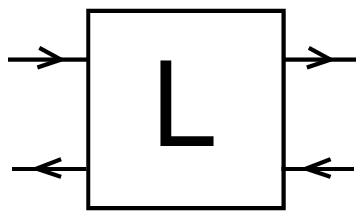}\\
(a) & (b) & (c)\\
\end{tabular}
\par\vspace{0.5cm}
Fig. 0.5
\end{center}

We perform a mutation of the link, of which $L$ is a part, 
by rotating the tangle along the 
$x, y$ or $z$ coordinate axis by the angle $\pi$. Thus we have three 
mutations $m_x, m_y$ and $m_z$, respectively; Fig. 0.6. 
Notice that together with the identity map they form the group 
$D_2=Z_2\oplus Z_2$. We keep the part of the link outside the tangle 
fixed and, if necessary change the orientation of the tangle part so it 
agrees with outside part of the link.

\par\vspace{1cm}
\begin{center}
\begin{tabular}{c}
\includegraphics[trim=0mm 0mm 0mm 0mm, width=.35\linewidth]
{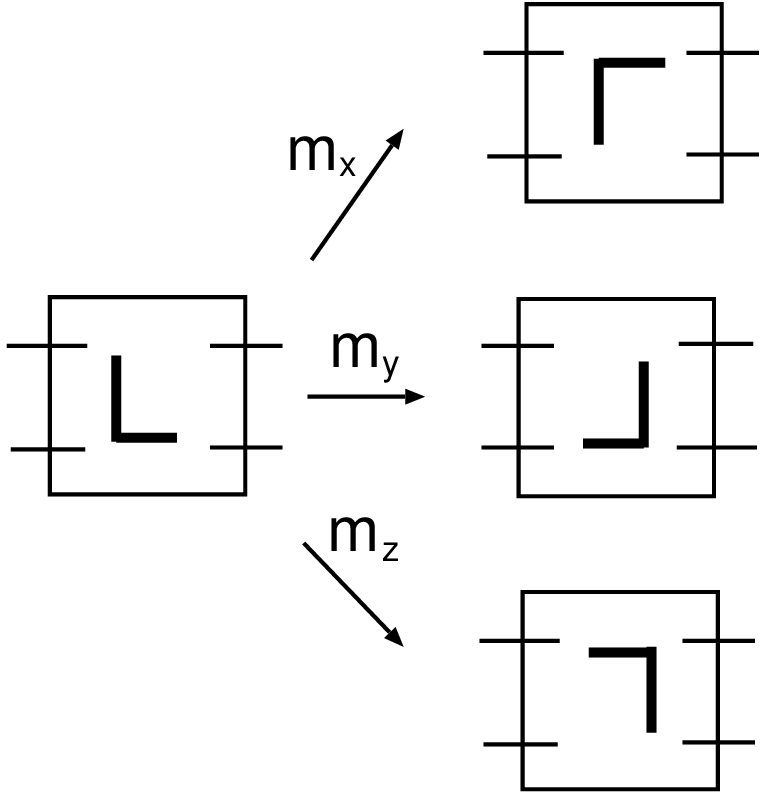}\\
\end{tabular}
\par\vspace{0.5cm}
Fig. 0.6
\end{center}

The mutation preserves not only Jones 
type polynomials but also the volume of the (hyperbolic) complement of
a link and 
the homeomorphic type of the branched double cover of $S^3$ with the
link as the branching set, among many other invariants.
The simplest pair of non-equivalent mutant knots, the Conway and
Kinoshita-Terasaka knots, is drawn in Fig. 0.2.

To show that Jones
type invariants are preserved by mutation, one follows Conway's  idea of skein
theory \cite{Co,L-M,Hos}.  Namely,  one simplifies the 2-tangle of
mutation as far as possible reaching, finally, tangles which are invariant
under mutation. In the case of the Kauffman skein relation one gets one of
the tangles of Fig. 0.7 with possible additional trivial components. 
Of course, each of these tangles is invariant under any mutation. 

\par\vspace{1cm}
\begin{center}
\begin{tabular}{ccc}
\includegraphics[trim=0mm 0mm 0mm 0mm, width=.2\linewidth]
{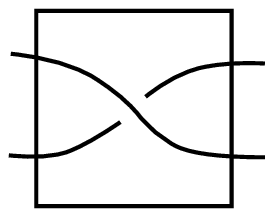}
&
\includegraphics[trim=0mm 0mm 0mm 0mm, width=.18\linewidth]
{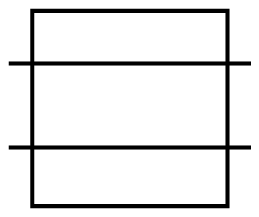}
&
\includegraphics[trim=0mm 0mm 0mm 0mm, width=.2\linewidth]
{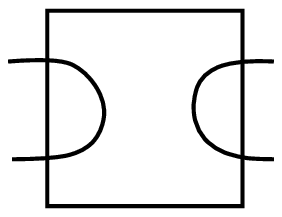}\\
\end{tabular}
\par\vspace{0.5cm}
Fig. 0.7
\end{center}

\section{Satellites of mutants}\label{1}

One can produce more complicated links from the given one, say $L$, by 
decorating each component of the link by some pattern. The resulting link 
is called a {\it satellite} of $L$.
If we consider satellites of mutants (with the same pattern) we obtain 
links with the same Jones polynomial and, sometimes, with the same skein
and Kauffman polynomials.\\
We should stress, however, that we cannot produce, in such a way, 
a nontrivial knot with a polynomial of the trivial knot. It is the case
because a mutation of a trivial knot is a trivial knot (the 2-fold 
branched cover of $(S^3,L)$ and $(S^3,m(L))$ are homeomorphic \cite{
Mon,Viro}, and
the 2-fold branched cover of $S^3$ with a nontrivial knot as the
branching set cannot be $S^3$; \cite{Wa}). Furthermore a nontrivial 
satellite of a nontrivial knot is a nontrivial knot.
\begin{theorem}[\cite{M-T-1}]\label{1.1} 
Let $D$ be a diagram of a framed unoriented link and let $m(D)$ be a 
mutant of $D$. Assume additionally that the mutation preserves
link components of $D$ (i.e. if $v$ is a boundary point of the 
rotated tangle, then $v$ and $m(v)$
lie in the same link component of $D$). 
Then for any satellites of $D$ and $m(D)$ with the same pattern 
$s$ (denoted by $s(D)$ and $s(m(D))$, respectively), the Kauffman
bracket polynomial is the same. 
That is $<s(D)>=<s(m(D))>$. If $D$ is additionally
oriented then $V_{s(D)}(t)=V_{s(m(D))}(t)$.
\end{theorem}

\begin{theorem}[\cite{L-L,P-1}]\label{1.2}
Let $m(D)$ be a mutant of $D$ obtained by rotating a tangle $T$ and
preserving link components of $D$. 
Let $s$ be a satellite pattern with
the wrapping number at most 2 (e.g. 2-cable or Whitehead double). Then:
\begin{enumerate}
\item
[(a)] If $D$ is oriented then:
\begin {enumerate} 
\item
[(i)] $P_{s(D)}(v,z) = P_{s(m(D))}(v,z)$, 
\item
[(ii)] $F_{s(D)}(a,x) = F_{s(m(D))}(a,x)$.
\end{enumerate}
\item
[(b)] If $D$ is an unoriented framed diagram, then\\
${\Lambda}_{s(D)}(a,x) = {\Lambda}_{s(m(D))}(a,x)$.
\end{enumerate}
\end{theorem}

Theorem 1.2 does not hold for wrapping number greater than 2. For example,
Morton and Traczyk  \cite{M-T-1} and J. Murakami \cite{Mura} 
have shown that a 3-cable of the Conway knot and its mutant 
the Kinoshita-Terasaka knot have different skein
and Kauffman polynomials. However a weaker fact still holds:
\begin{theorem}[\cite{L-bull,P-1}]\label{1.3}\
\\
Let $L_1$ and $L_2$ be two oriented links and  let $-L_2$ be obtained from 
$L_2$ by reversing  orientations of all components of $L_2$.  
Consider connected sums
$L_1\# L_2$ and $L_1\# -L_2$, where the same components of $L_1$ and 
$L_2$ are involved in $L_1\# L_2$ and $L_1\# -L_2$ (the second sum 
can be thought as a degenerated, components preserving, mutation of the first). Then for any pattern $s$, the satellites
$s(L_1\# L_2)$ and $s(L_1\# -L_2)$ have the same skein and Kauffman
polynomials.
\end{theorem}

The original proofs of Theorems 1.1-1.3 were combinatorial, using 
skein theory similarly as in the proof of the case of mutation. Later,
however, J. Murakami and G.Kuperberg found proofs based on 
properties of irreducible representations of Lie algebras.

Theorems 1.1 and 1.3 are the basic tools in constructing different
3-manifolds with the same Witten invariants \cite{Kania,Li-1}.
Theorem 1.1 is used in the case of $SU(2)$ Witten invariant, as
constructed by Reshetikhin and Turaev \cite{R-T} and Theorem 1.3
in the case of classical Lie algebra Witten invariants, as constructed
by Turaev and Wenzl \cite{T-W}.

\section{Rotors}\label{2}
 The idea of rotors was used first in graph theory in the fundamental
paper on division of a square into smaller unequal squares \cite{BSST}.
Later Tutte used it to produce different graphs with the same 
dichromatic polynomial \cite{Tut}. The idea was ``translated" into knot 
theory in the paper of R.Anstee, D.Rolfsen and myself.
\begin{definition}[\cite{A-P-R}]\label{2.1}\
\\
Consider an $n$-tangle, that is a part of the link diagram (possibly
oriented or framed - depending on the application),  placed in the regular
$n$-gon with $2n$ boundary points (n inputs and n outputs). We say that
this $n$-tangle is an $n$-rotor if it has a rotational symmetry, that is
the tangle is invariant with respect to rotation along $z$-axis by
the angle $\frac{2\pi}{n}$; see Fig. 2.1 for an example of a 4-rotor.
\end{definition}

\par\vspace{1cm}
\begin{center}
\begin{tabular}{c}
\includegraphics[trim=0mm 0mm 0mm 0mm, width=.35\linewidth]
{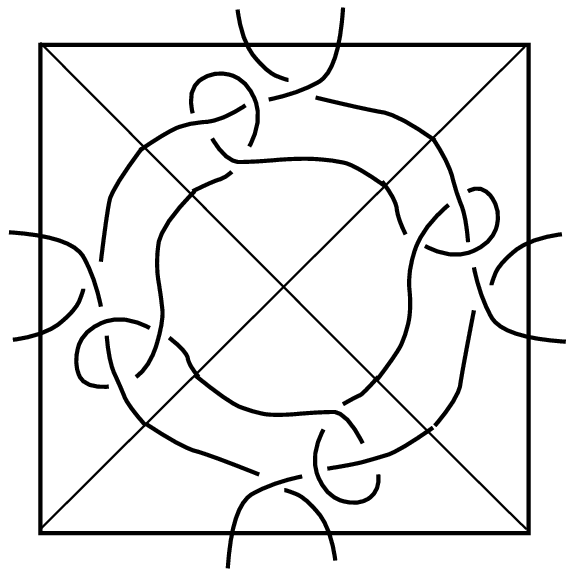}\\
\end{tabular}
\par\vspace{0.5cm}
Fig. 2.1. Rotor of Toulouse, Carolingian era, IX century AD.
\end{center}

\begin{theorem}[\cite{A-P-R}]\label{2.2} \ \\
Let $L$ be a link diagram with an n-rotor part $R$. Let the rotant $r(L)$, 
be obtained from $L$ by rotating $R$ along the $x$-axis by the angle $\pi$
and keeping the stator, $S=L-R$, unchanged (if necessary, we change the 
orientation of the rotor so it agrees with that of the stator). Then:
\begin{enumerate}
\item
[(a)] $<L>=<r(L)>$ for $n\leq 5$, where $L$ is an unoriented framed 
link diagram.
\item
[(b)] $P_L(v,z)=P_{r(L)}(v,z)$ for $n\leq 4$, where  $L$ is an 
oriented link diagram.
\item
[(c)]  ${\Lambda}_L(v,x)={\Lambda}_{r(L)}(v,x)$ for $n\leq 3$, where 
$L$ is an unoriented framed link diagram. 
\end{enumerate}
Furthermore, if $L$ is oriented then $Tait(L)=Tait(r(L))$.
\end{theorem}

The proof of the theorem given in \cite{A-P-R} is a straightforward 
generalization of the proof for mutants: the skein relation allows as to
simplify the stator of $L$ so that it is invariant under the reflection
in a side of the $n$-gon. We illustrate it for $n=3$ and the Kauffman skein
relation. If we simplify the stator (placed in a regular triangle) using
the Kauffman skein relations, we obtain one of the fifteen 3-tangles 
(possibly with additional trivial components); Fig. 2.2.

\par\vspace{1cm}
\begin{center}
\begin{tabular}{c}
\includegraphics[trim=0mm 0mm 0mm 0mm, width=.8\linewidth]
{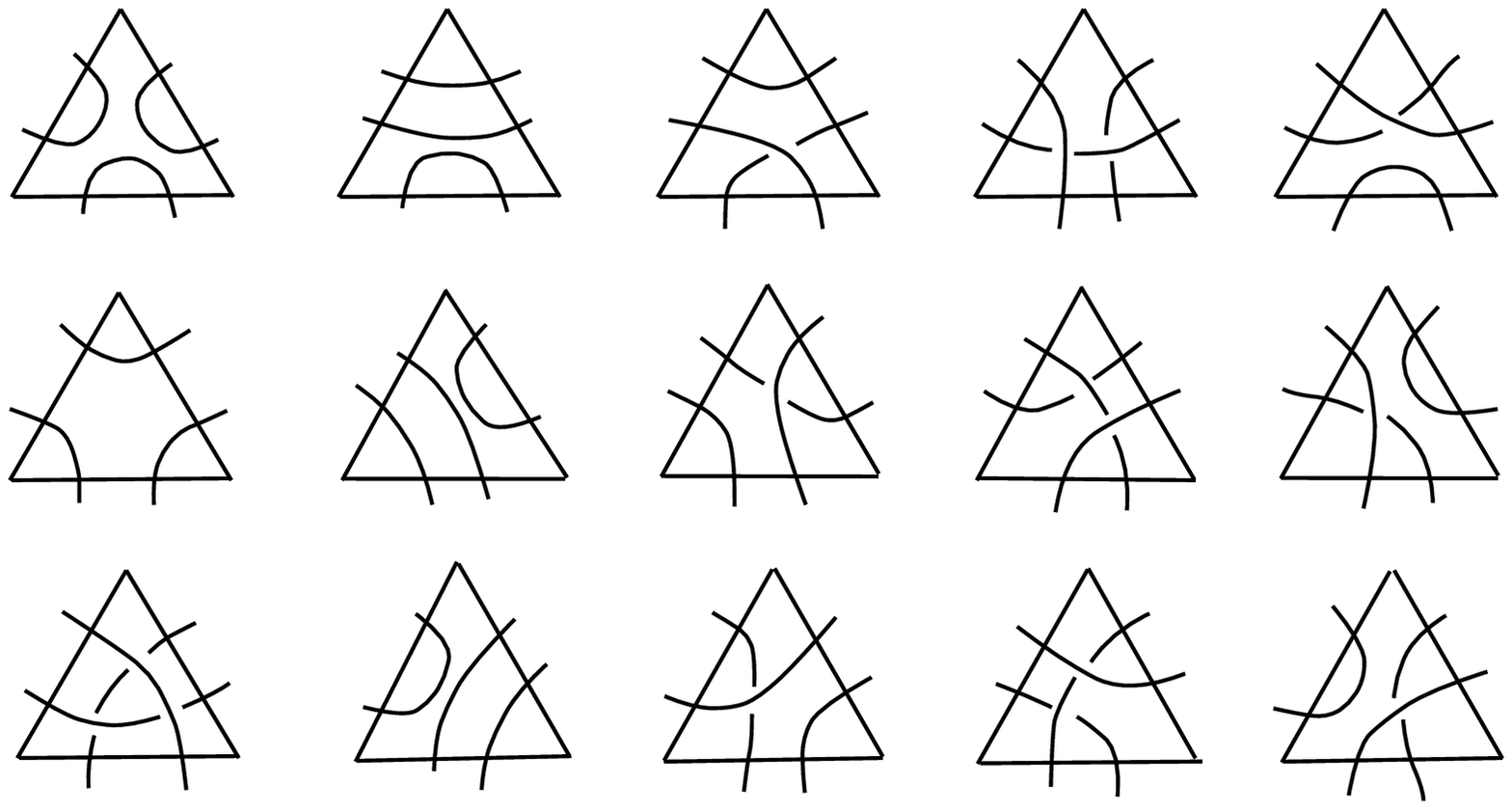}
\\
\end{tabular}
\par\vspace{0.5cm}
Fig. 2.2
\end{center}

Clearly each tangle of Fig. 2.2 is symmetric with respect to reflection
in an edge of the triangle\footnote{Added for e-print: Equivalently 
for every 3-tangle in Fig. 2.2 there is an axis of symmetry of the 
triangle such that the tangle is invariant under rotation along this axis 
by the angle $\pi$.}. Instead of reflecting the stator we can 
reflect the rotor (and because of its rotational symmetry, all reflections
are equivalent), and the Theorem 2.2(c) follows.

D.Rolfsen was searching for examples of $n$-rotors for which the theorem
does not hold. In particular he studied 6-rotors and their Jones polynomial.
However he couldn't find a counterexample to the theorem. At the Sussex
conference held in summer of 1987, he discussed rotors with P.Traczyk and
Traczyk observed that 6-rotors considered by Rolfsen have only
2 connections between identical segments of rotors; see Fig. 2.3.
\\

\begin{center}
\begin{tabular}{c}
\includegraphics[trim=0mm 0mm 0mm 0mm, width=.35\linewidth]
{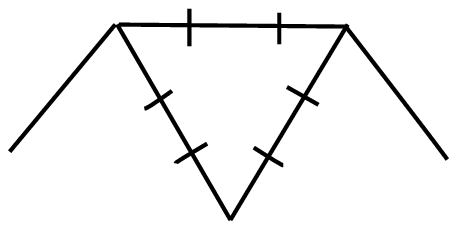}\\
\end{tabular}
\par\vspace{0.5cm}
Fig. 2.3
\end{center}

T.Jin and Rolfsen considered more complicated rotors and  found, in 
the summer
of 1988, that the Theorem 2.2 is the best possible \cite{J-R}.\\
On the other hand, Traczyk showed that if there are only two internal
connections between segments of a rotor then the Jones and
the skein polynomials are the same for the link and its rotant (he proved
it first for $n=6$, in 1987, and the next year in the full generality). 
Namely: 
\begin{theorem}[\cite{Tr-1}]\label{2.3}\
\\
Let $L$ be a link diagram with an n-rotor part $R$. Furthermore assume that
there are at most two arcs between neighboring segments of the rotor.
If $r(L)$ is the rotant of $L$ then:
\begin{enumerate}
\item
[(a)] $<L>=<r(L)>$, for any $n$, where  $L$ is an unoriented framed link 
diagram. 
\item
[(b)] $P_L(v,z)=P_{r(L)}(v,z)$, for any $n$, where $L$ is an oriented link 
and a segment of the rotor (which is a 3-tangle) is oriented, 
up to the global change of orientation, as in Fig. 2.4.
\end{enumerate}
\end{theorem}

\par\vspace{1cm}
\begin{center}
\begin{tabular}{c}
\includegraphics[trim=0mm 0mm 0mm 0mm, width=.35\linewidth]
{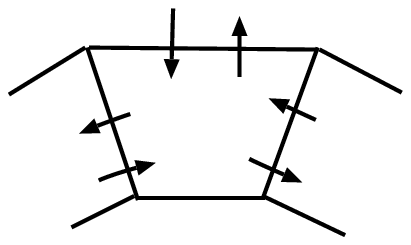}\\
\end{tabular}
\par\vspace{0.5cm}
Fig. 2.4
\end{center}

The method used by Traczyk for proving Theorem 2.3, is essentially different
than that used in the proof of Theorem 2.2. Namely, it is algebraic and 
uses essentially the linearity of skein relations. Furthermore Traczyk 
operates on the rotor part of the link (excluding its center)
instead  of the stator, so he gets, in fact, theorem about tangles in the
solid torus (projected into an annulus), or any 3-manifold in which the
solid torus is embedded (to make it precise one should consider the notion
of the skein module of a 3-manifold; see part 3).
Traczyk's theorem does not hold for the Kauffman polynomial. A simple
counterexample is given in \cite{J-R} (Example 1); see Fig. 2.5, where
4-rotors with different Kauffman polynomials are presented. Furthermore
these rotors have only two connections between segments.

\par\vspace{1cm}
\begin{center}
\begin{tabular}{cc}
\includegraphics[trim=0mm 0mm 0mm 0mm, width=.35\linewidth]
{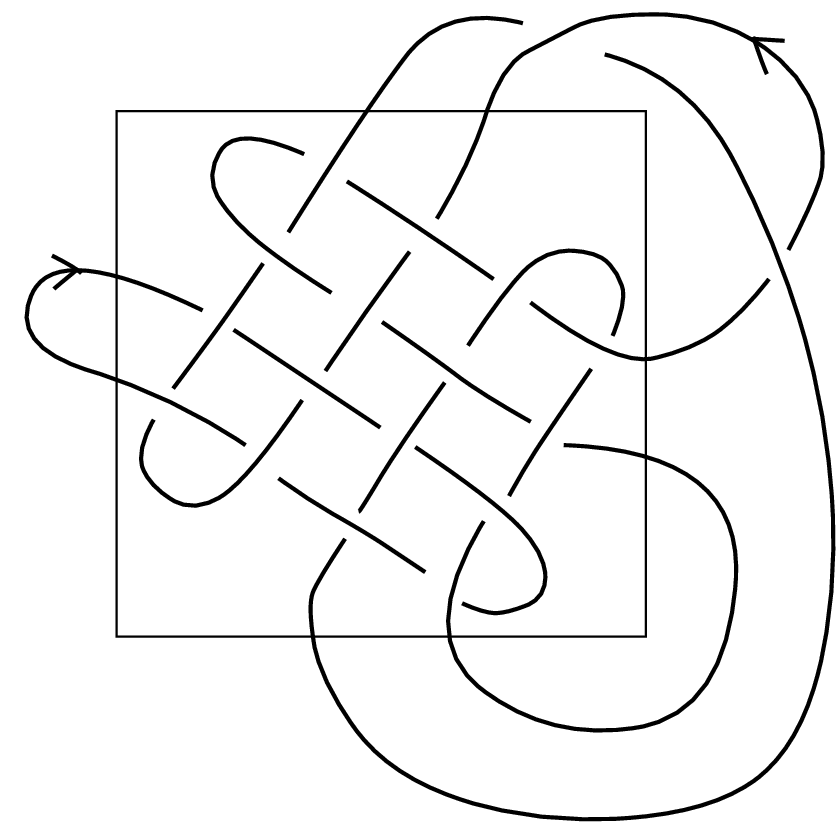}&
\includegraphics[trim=0mm 0mm 0mm 0mm, width=.38\linewidth]
{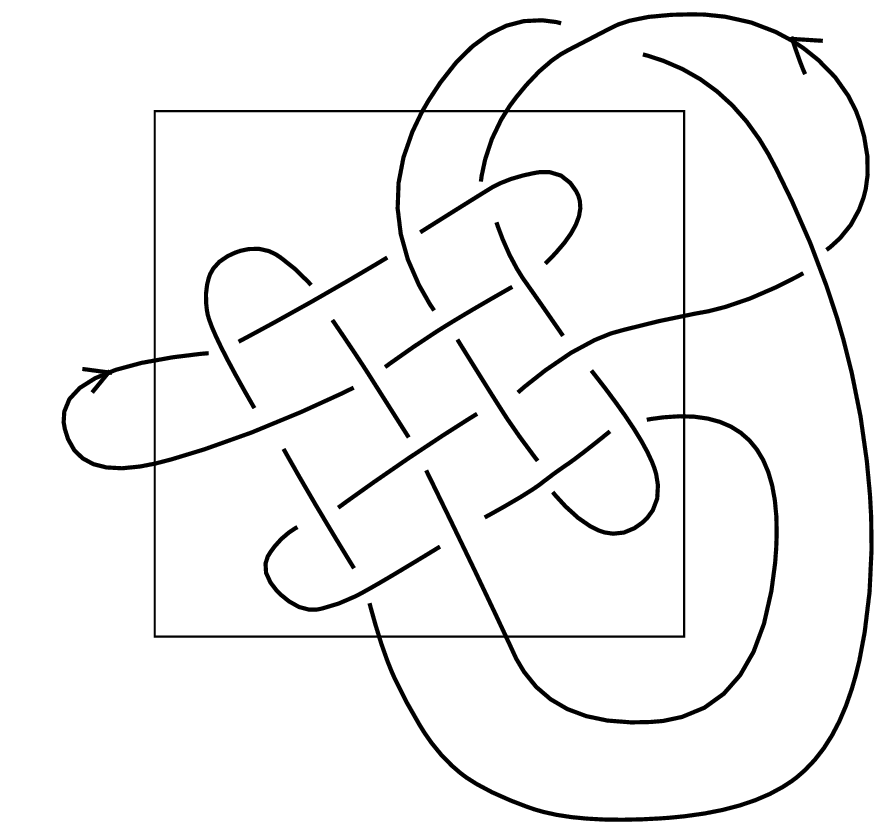}
\\
\end{tabular}
\par\vspace{0.5cm}
Fig. 2.5. $F_L(a,x)-F_{r(L)}(a,x)=(a^{10}+4a^{12}+6a^{14}+4a^{16}+
a^{18})+ {\cal O}(z).$
\end{center}
\ \\

It follows from the work of Tutte that the determinant of an 
oriented link, $|{\Delta}_L(-1)|$, is the same for $L$ and its 
rotant $r(L)$, for any $n$ 
and any number of connections between segments. This suggests 
the possibility that the whole Alexander polynomial is unchanged by
rotation and a lot of examples which confirm this were computed
by Jin and Rolfsen.
\begin{conjecture}\label{2.4}
If $L$ is an oriented link and $r(L)$ its rotant then their Alexander
polynomials are equal\footnote{Added for e-print: The Jin-Rolfsen 
conjecture has been proved by P.Traczyk for $n$-rotors with 
inputs and outputs alternating (Conway polynomial and oriented 
rotant links, {\it Geometria Dedicata}, 2004, to appear. Makiko Ishiwata 
constructed a counterexample using a 6-rotor with ``2 inputs, 2 outputs..." 
pattern (Rotation and Link Invariants, Doctorate Thesis, OCU, Japan, 2004, 
and M.D{\c a}bkowski, M.Ishiwata, J.H.Przytycki, A.Yasuhara, 
Rotation and signature invariants, preprint 2003).}. 
\end{conjecture}
We noticed in \cite{A-P-R}, that if we decorate an unoriented framed
 link $L$ and its 3-rotant $r(L)$ by a pattern $s$ with the
 wrapping number 2, 
then the Kauffman bracket polynomials of $s(L)$ and $s(r(L))$ are equal.
There is no need, however, for a separate proof of this fact 
because Yamada showed
\cite{Ya} that the Kauffman bracket of a 2-cable of a link is determined 
by the Kauffman polynomial of the link at $(a,x)=(iA^8,-i(A^4-A^{-4})).$ 
The idea of rotors was taken from graph theory\footnote{
One can say that it came from physics, as the authors of \cite{BSST} 
where motivated by the theory of electrical circuits.}. 
Now however, Traczyk's theorem and the results of the next parts 
can be translated back to give new results in graph theory.

\section{The spectral parameter tangle}\label{3}
 In the spring of 1991, Jones used an idea from the statistical mechanics 
to produce examples of different links with the same Jones type
polynomials. 
We describe in this part of the talk, 
how the Yang-Baxter equation with
spectral parameter can be ``translated" into an equation involving tangles
and how the Baxter method of "commuting transfer matrices"
can be "translated" to produce various links with the same Jones,
and skein polynomials \cite{Jo-conf,Jo-3}. 
Subsequently we generalize slightly the method of Jones, and relate it to 
Traczyk's work on rotors of links.\footnote{This part and Part 4 are
based on the notes for the talk given at the University
of Tennessee, October 1991 \cite{P-3}.} Finally we give another application
of the spectral parameter tangles.

   We first present the work of Jones.
 
Consider a finite dimensional vector space $V$ and the space
of parameters $\Lambda$. Consider also the family, $R(\lambda)$, of
matrices $(R: \Lambda \to End(V\otimes V))$. We say that $R$ satisfies the 
Yang Baxter equation with spectral parameter if for every pair
$\lambda, {\lambda}'\in \Lambda$ there is $\lambda'' \in \Lambda$ such 
that $$R_1(\lambda)R_2(\lambda ')R_1(\lambda '') = 
R_2(\lambda '')R_1(\lambda ')R_2(\lambda),$$ 
where both sides are endomorphisms of $V\otimes V\otimes V$
and $R_1(\chi) = R(\chi)\otimes Id$  and $R_2(\chi)=Id\otimes R(\chi).$\ \
$\chi \in \Lambda$ is called the {\em spectral
parameter}. The Yang Baxter equation is expressed graphically in Figure 3.1.

\par\vspace{1cm}
\begin{center}
\begin{tabular}{c}
\includegraphics[trim=0mm 0mm 0mm 0mm, width=.8\linewidth]
{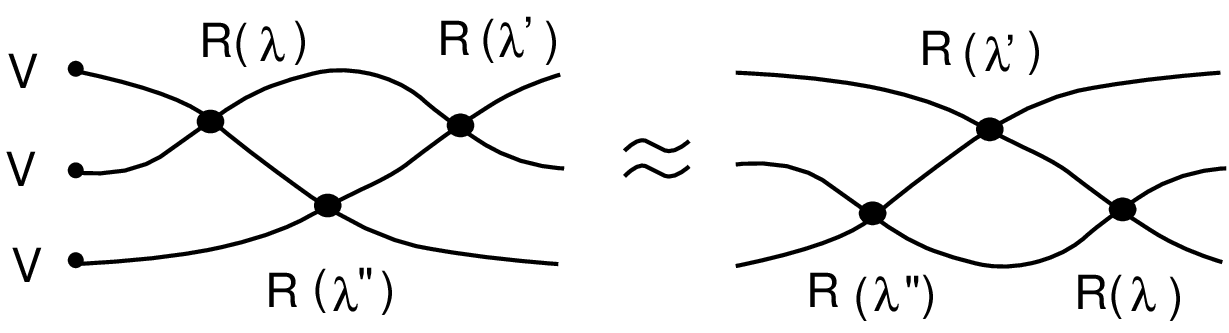}\\
\end{tabular}
\par\vspace{0.5cm}
Fig. 3.1
\end{center}

We interpret three initial points as $V\otimes V\otimes V$ and
going through the crossing corresponds to the endomorphism $R(\chi)$. 
When the spectral parameter is constant then the correspondence 
between the crossing \ and the
Yang-Baxter operator $R$ can be used to define new, Jones type, invariants
of links (Jones \cite{Jo-2}, Turaev \cite{Tu}). 
The method uses the fact that a
Yang-Baxter equation corresponds to the third Reidemeister
move (Fig. 3.2).

\par\vspace{1cm}
\begin{center}
\begin{tabular}{c}
\includegraphics[trim=0mm 0mm 0mm 0mm, width=.7\linewidth]
{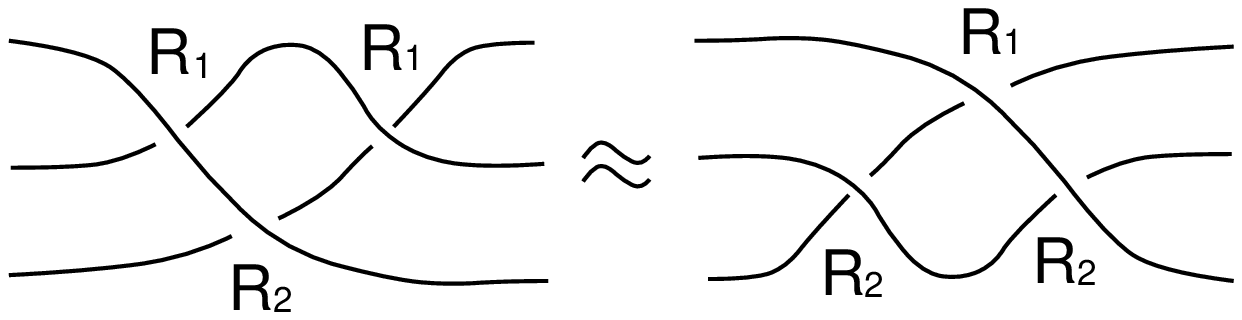}\\
\end{tabular}
\par\vspace{0.5cm}
Fig. 3.2
\end{center}
 
V.Jones  in the spring of 1991 \cite{Jo-conf,Jo-3} discovered 
that one can use the idea of the Yang-Baxter operator with spectral  
parameter to produce different links with the same Jones 
polynomials.
The idea uses the classical argument from the theory of solvable
models in statistical mechanics \cite{Baxter}. As before, $V$ 
corresponds to a point
and $\underbrace{V\otimes V\otimes ...\otimes V}_n$ to $n$ points. 
$R(\lambda) \in
End(V\otimes V)$  corresponds to a 2-tangle  $L$   and the composition of 
endomorphisms corresponds to the composition of tangles. 
To get the analogy to the Yang-Baxter equation we consider the skein 
module (linear skein) of a tangle,  
the notion which was essentially introduced by J.Conway 
(\cite{Co,P-2,Tu-1,H-P-1}).
In short skein modules are quotients of modules over ambient
isotopy classes of links in a 3-manifold (possibly with boundary)
 by properly chosen local (skein) relations.

Now we have to be more specific and choose a skein module with which we
work. We start with the Kauffman bracket approach to the Jones 
polynomial. 
Consider an $n$-tangle, that is a 2-disc with fixed $2n$ points on the 
boundary and a  framed link diagram (composed of closed curves and 
curves with endpoints fixed) inside. Framing is fixed at boundary points. 
The Kauffman bracket skein module of n-tangles, ${\cal S}_{2,\infty}(n)$
is defined to be an $R$-module (for a chosen commutative ring with 1) obtained
from the free $R$-module over all tangles up to isotopy (modulo boundary)
divided by the submodule generated by (the Kauffman bracket) skein relations 
$L_A = AL_0 + A^{-1}L_{\infty}$ \ and framing relations
 $L^{(1)} = -A^3L $, where $L^{(1)}$ denotes a framed link obtained from $L$
by twisting its framing once in the positive direction. \
 $A$ is an invertible element of $R$ (in practice, unless otherwise stated, 
we will assume at this talk that $R= {\cal F}(A)$ - the field of rational 
functions in variable $A$).
${\cal S}_{2,\infty}(n)$ is a free module of 
$$\frac{1}{n+1}{{2n}\choose{n}}$$ 
generators (\cite{Jo-1,M-T-2}).\ \  
For example ${\cal S}_{2,\infty}(2)$ has two 
generators  and ${\cal S}_{2,\infty}(3)$ has five generators; Fig. 3.3.

\par\vspace{0.5cm}
\begin{center}
\begin{tabular}{c}
\includegraphics[trim=0mm 0mm 0mm 0mm, width=.85\linewidth]
{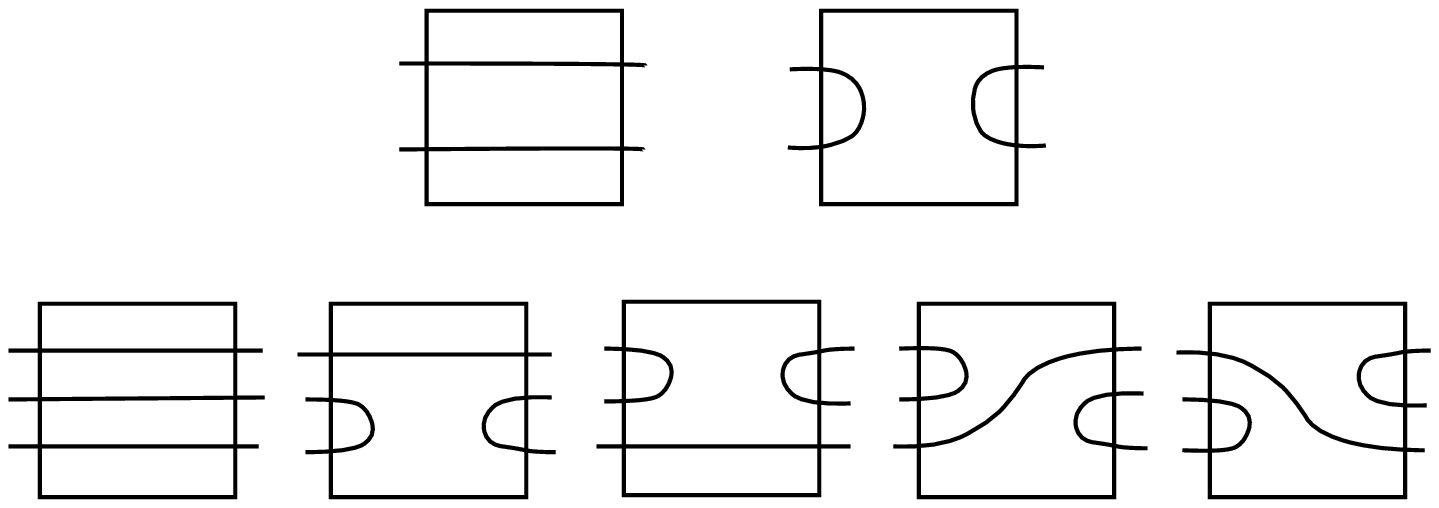}\\
\end{tabular}
\par\vspace{0.2cm}
Fig. 3.3
\end{center}

${\cal S}_{2,\infty}(n)$ with product yielded be 
the composition of tangles has a structure of an algebra 
(the Temperley-Lieb algebra \cite{T-L}).\ 
We will use the standard (in Yang-Baxter equations theory) notation:\
$R_{(i)}$ for a 2-tangle $R$ placed on $i$'th and $(i+1)$'th strings of 
$n$ strings, see Fig. 3.4.
The result (analogy to the Yang-Baxter equation with spectral 
parameter) which Jones uses, can be stated as follows:
\begin{lemma}[Jones \cite{Jo-3}]\label{Lemma 3.1}
\ \\
For a dense subset of pairs of tangles 
$T,T' \in {\cal S}_{2,\infty}(2)\times
{\cal S}_{2,\infty }(2)$ there is an invertible tangle $T''$
(in Temperley-Lieb algebra ${\cal S}_{2,\infty}(2)$) such that:
$T_{(1)}T'_{(2)}T''_{(1)} = T''_{(2)}T'_{(1)}T_{(2)}$.
Graphically this equality of two 3-tangles is shown in Figure 3.4. 
\end{lemma}
\par\vspace{0.5cm}
\begin{center}
\begin{tabular}{c}
\includegraphics[trim=0mm 0mm 0mm 0mm, width=.8\linewidth]
{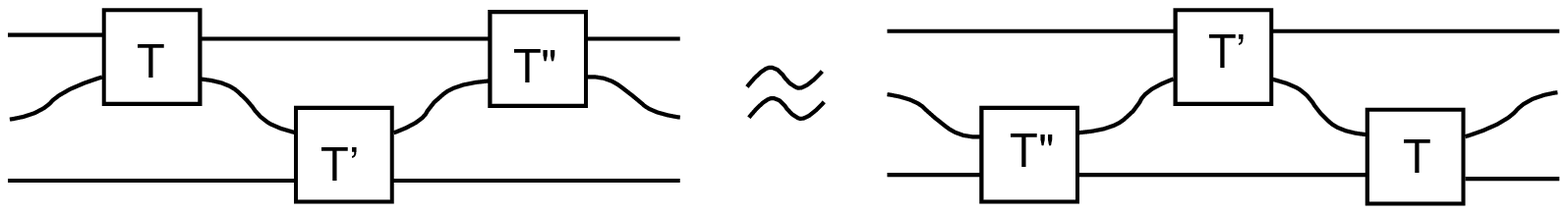}\\
or equivalently
\\
\includegraphics[trim=0mm 0mm 0mm 0mm, width=.6\linewidth]
{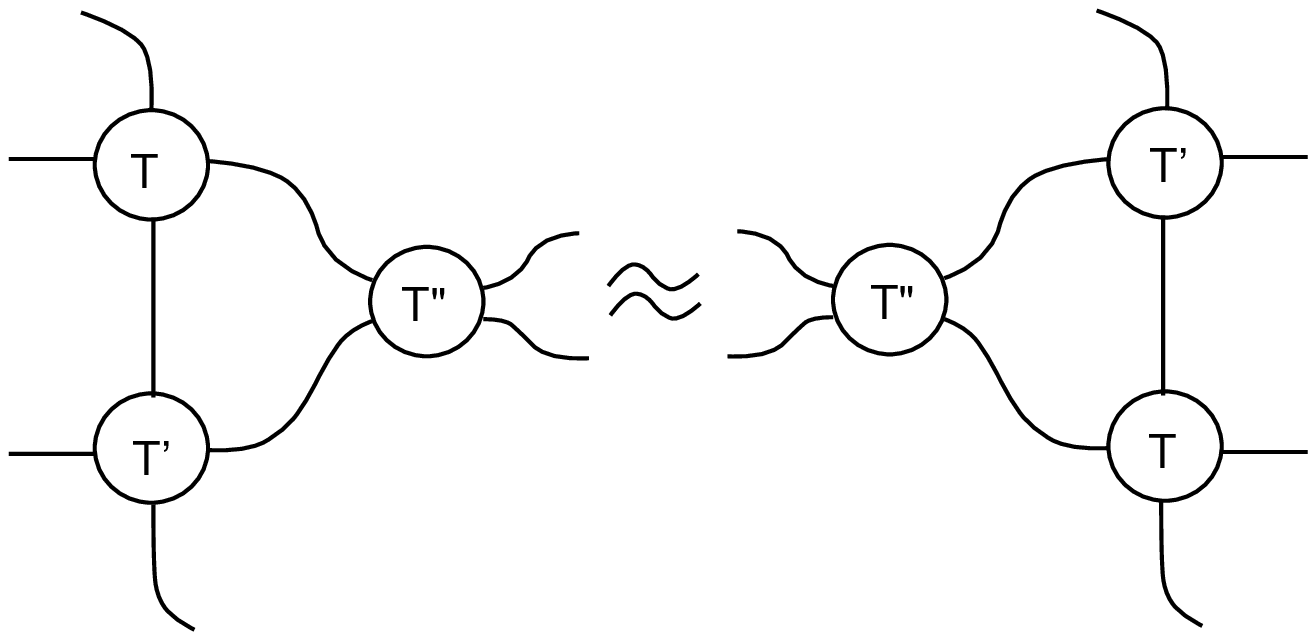}
\\
\end{tabular}
\par\vspace{0.2cm}
Fig. 3.4
\end{center}

We do not prove here the Jones lemma as its generalized version
is given in Lemma 3.3. 
The ``trick", which Jones uses, goes as follows:\\
Consider 2-tangles $T$ and $T'$ placed cyclicly in the annulus
as in Figure 3.5(a). Assume that for $T$ and $T'$, there exists
an invertible 2-tangle $T''$ from Lemma 3.1. Then without
changing the element of the Kauffman bracket skein module of the 
solid torus (annulus is the projection surface of the solid torus),
we can insert the pair of tangles $T''$ and $(T'')^{-1}$ (Fig. 3.5(b))
and then move around the annulus with $T''$, interchanging $T$ with
$T'$ on the way (Fig. 3.5(c)), and finally arriving on the second
side of $(T'')^{-1}$ and canceling it; Fig. 3.5(d).

\par\vspace{1cm}
\begin{center}
\begin{tabular}{c}
\includegraphics[trim=0mm 0mm 0mm 0mm, width=.7\linewidth]
{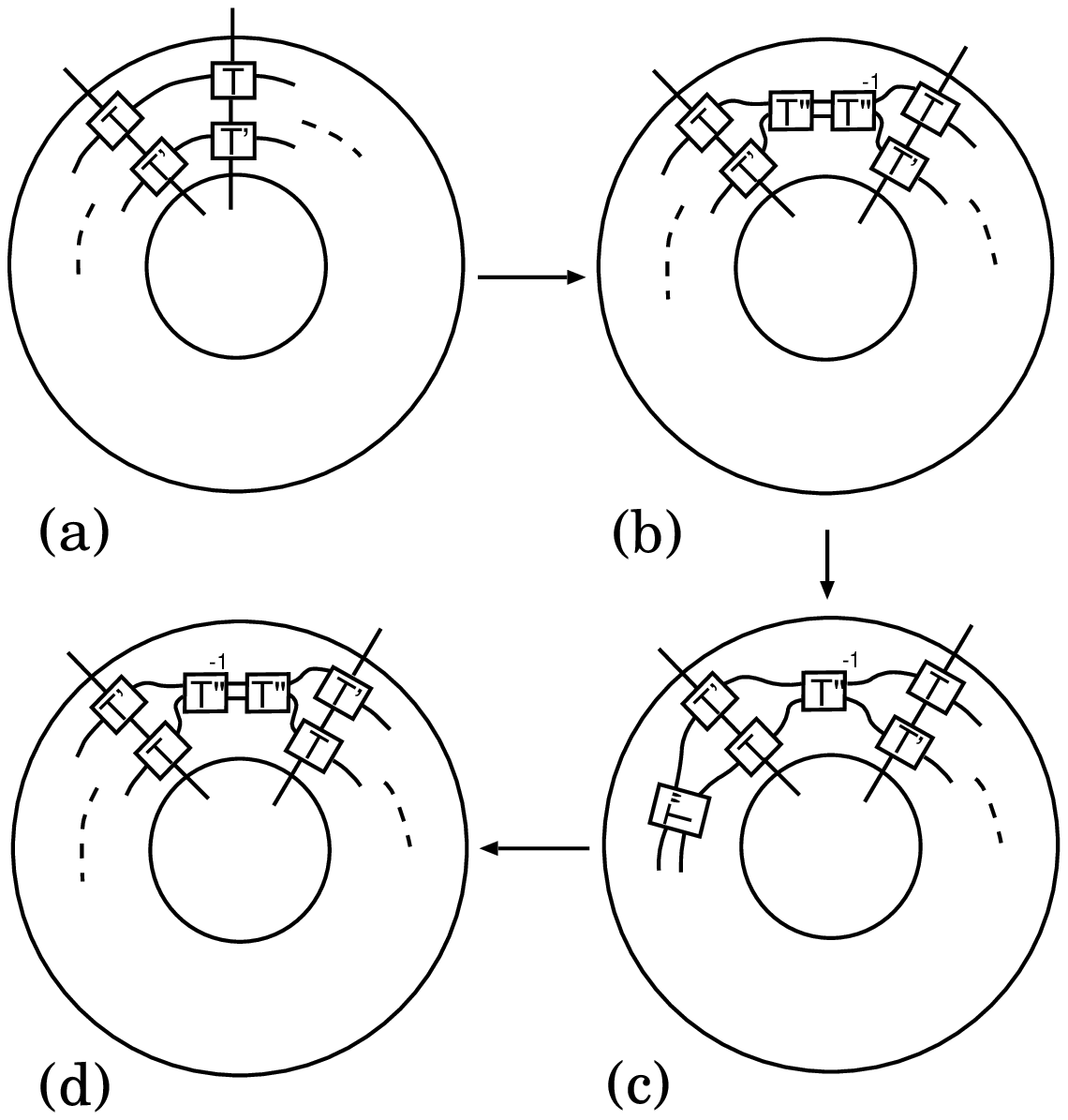}\\
\end{tabular}
\par\vspace{0.5cm}
Fig. 3.5
\end{center}

One slightly extend the scope of the Jones ``trick" by observing that 
${\cal S}_{2,\infty}(2)$ is commutative. Thus one gets.
\begin{theorem}[\cite{Jo-3}]\label{3.2}
Let $T$ and $T'$ be 3-tangles and $s_i$, i=1,...,n 2-tangles.
Consider the following elements of the Kauffman bracket skein module 
of the solid torus with 2n boundary points; Fig. 3.6(a) and (b).

\par\vspace{1cm}
\begin{center}
\begin{tabular}{c}
\includegraphics[trim=0mm 0mm 0mm 0mm, width=.9\linewidth]
{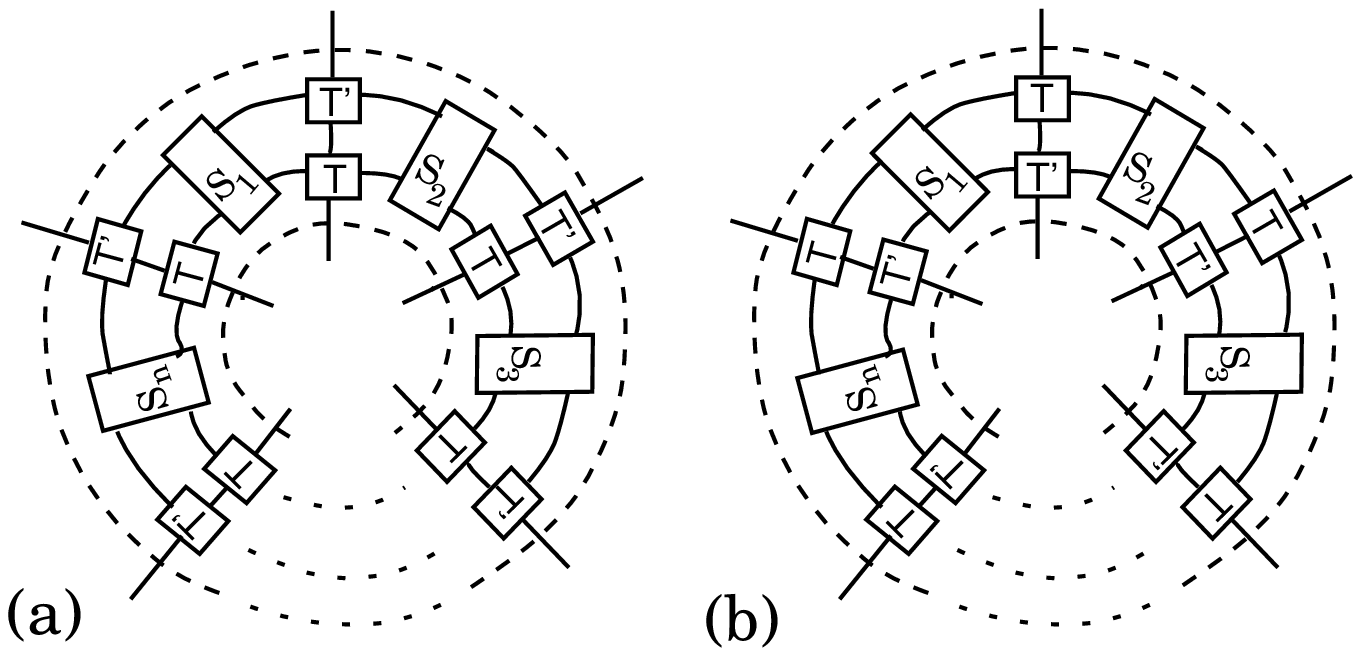}\\
\end{tabular}
\par\vspace{0.5cm}
Fig. 3.6
\end{center}

Then the tangles (a) and (b) represent the same element of the
Kauffman bracket skein module of the solid torus. In particular if
the solid torus is embedded in an oriented 3-manifold and boundary
points are connected together (outside the solid torus), then the links
obtained are equal in the Kauffman bracket skein module of the 3-manifold.
\end{theorem}

At the conference in Sacramento at April 1991, Jones gave a talk explaining
his ideas. We discussed them also afterwards (Jones, Hoste and myself).
I knew well the work of Traczyk on rotors so I suspected immediately that 
it should be related to results of Jones. In fact I noticed that the 
Jones method works without changes if instead of two tangles $T$ and 
$T'$ combined as in Figure 3.7(a), one can, more generally, consider one 
tangle of Figure 3.7(b). 
Similarly one can work with the tangle of Figure 3.7(c).
In such a way one recovers the Traczyk result and 
generalizes it. We present now, with details, 
this slight generalization of the Jones work.
We start from the case of the Kauffman bracket polynomial (and the
Kauffman bracket skein module).

\par\vspace{1cm}
\begin{center}
\begin{tabular}{c}
\includegraphics[trim=0mm 0mm 0mm 0mm, width=.9\linewidth]
{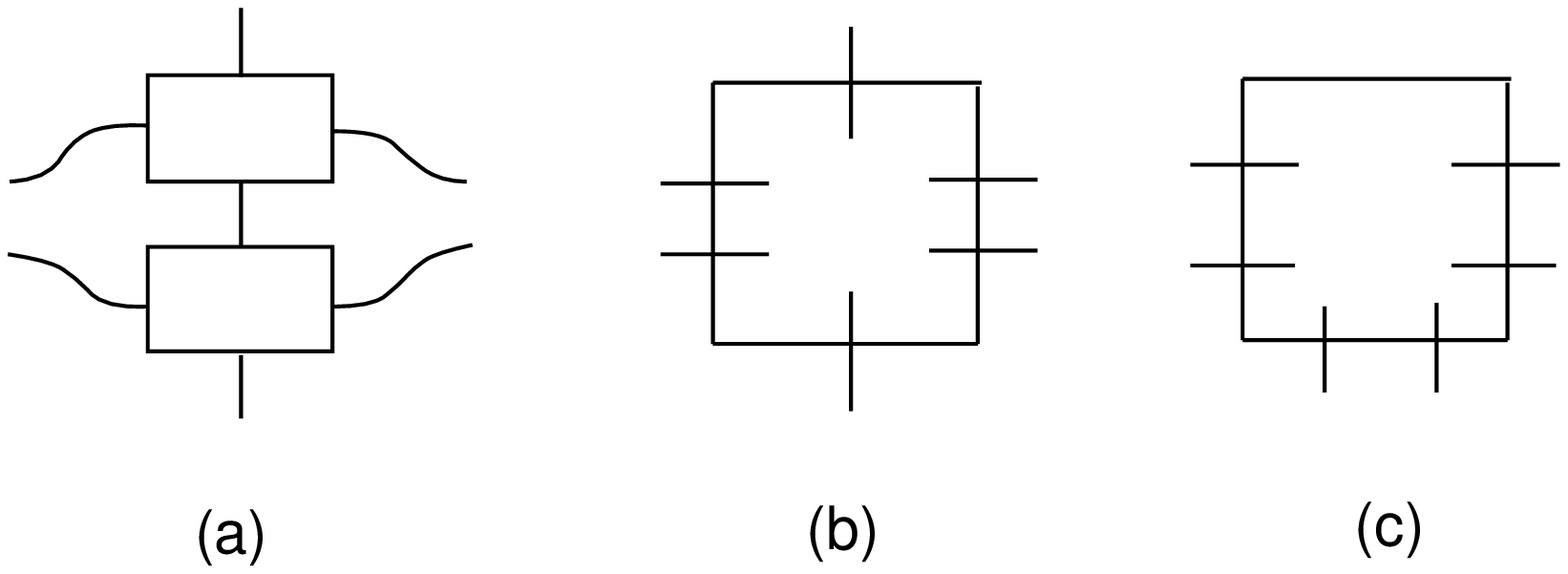}\\
\end{tabular}
\par\vspace{0.5cm}
Fig. 3.7
\end{center}

\begin{lemma}\label{Lemma 3.3}\
\begin{description}
\item
[(a)] 
Let $r_z : {\cal S}_{2,\infty}(n) \to {\cal S}_{2,\infty}(n)$ 
denotes the involution of modules (the algebra anti-isomorphism) 
 generated by rotating a tangle by 180 degrees about the $z$-axis\\
$(r_z(L) =\includegraphics[trim=0mm 0mm 0mm 0mm, width=.02\linewidth]
{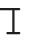}).$ Then for a dense subset of tangles 
$L \in {\cal S}_{2,\infty}(3)$
there is an invertible 2-tangle $P=P(L)$ such that
$LP_{(1)} = P_{(2)}r_z(L)$  in ${\cal S}_{2,\infty}(3)$  (see Figure 3.8).

\par\vspace{1cm}
\begin{center}
\begin{tabular}{c}
\includegraphics[trim=0mm 0mm 0mm 0mm, width=.8\linewidth]
{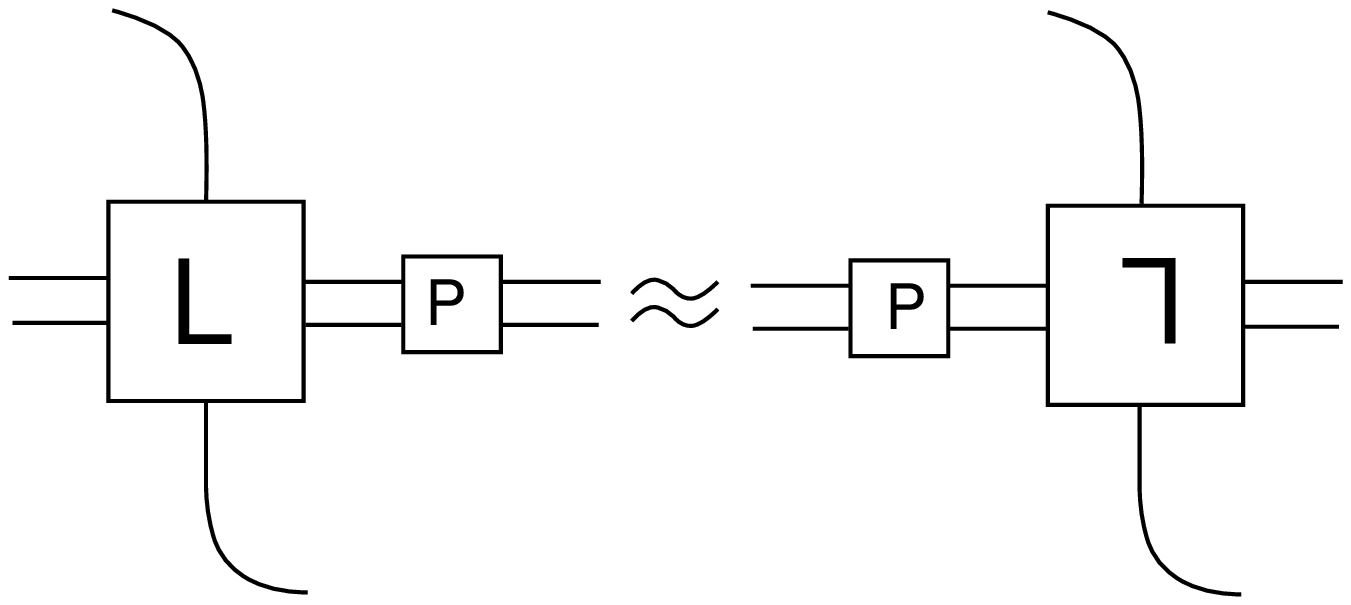}\\
\end{tabular}
\par\vspace{0.5cm}
Fig. 3.8
\end{center}

\item
[(b)]
Let $r_y: {\cal S}_{2,\infty}(n)\to {\cal S}_{2,\infty}(n)$ 
denote the involution of modules given
by the rotation $r_y$ ($r_y(L)=\includegraphics[trim=0mm 0mm 0mm 0mm, width=.013\linewidth]
{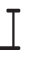}$). Then there is a dense subset, 
$D(3)$, of elements of ${\cal S}_{2,\infty}(3)$,
such that if $L$ is in this subset, then there is an invertible 2-tangle 
$P=P(L)$ such that $LP_{(1)} = P_{(1)}r_y(L)$  in ${\cal S}_{2,\infty}(3)$;  
(see Figure 3.9).
\end{description}
\end{lemma}

\par\vspace{1cm}
\begin{center}
\begin{tabular}{c}
\includegraphics[trim=0mm 0mm 0mm 0mm, width=.9\linewidth]
{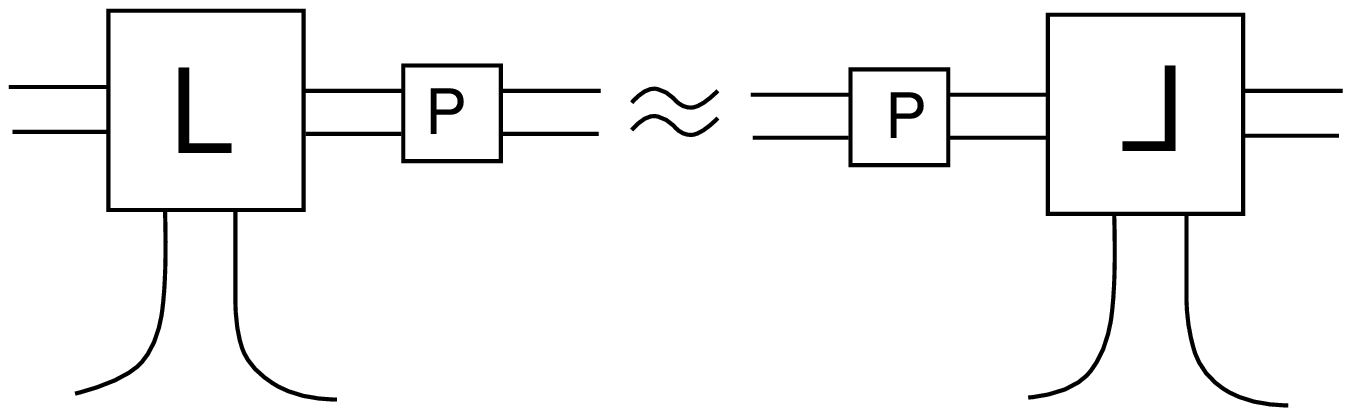}\\
\end{tabular}
\par\vspace{0.5cm}
Fig. 3.9
\end{center}

\begin{proof}
We will prove (b) in details. The proof of (a) is similar\footnote{
It is observed in \cite{H-P-2} that (a) follows immediately from (b),
thus there is no need for a separate calculation.}.
Let $(e_1,e_2,...,e_5)$ be the basis of 
${\cal S}_{2,\infty}(3)$ and $(f_1,f_2)$ 
the basis of ${\cal S}_{2,\infty}(2)$ as shown in Figure 3.10.

\par\vspace{1cm}
\begin{center}
\begin{tabular}{c}
\includegraphics[trim=0mm 0mm 0mm 0mm, width=.9\linewidth]
{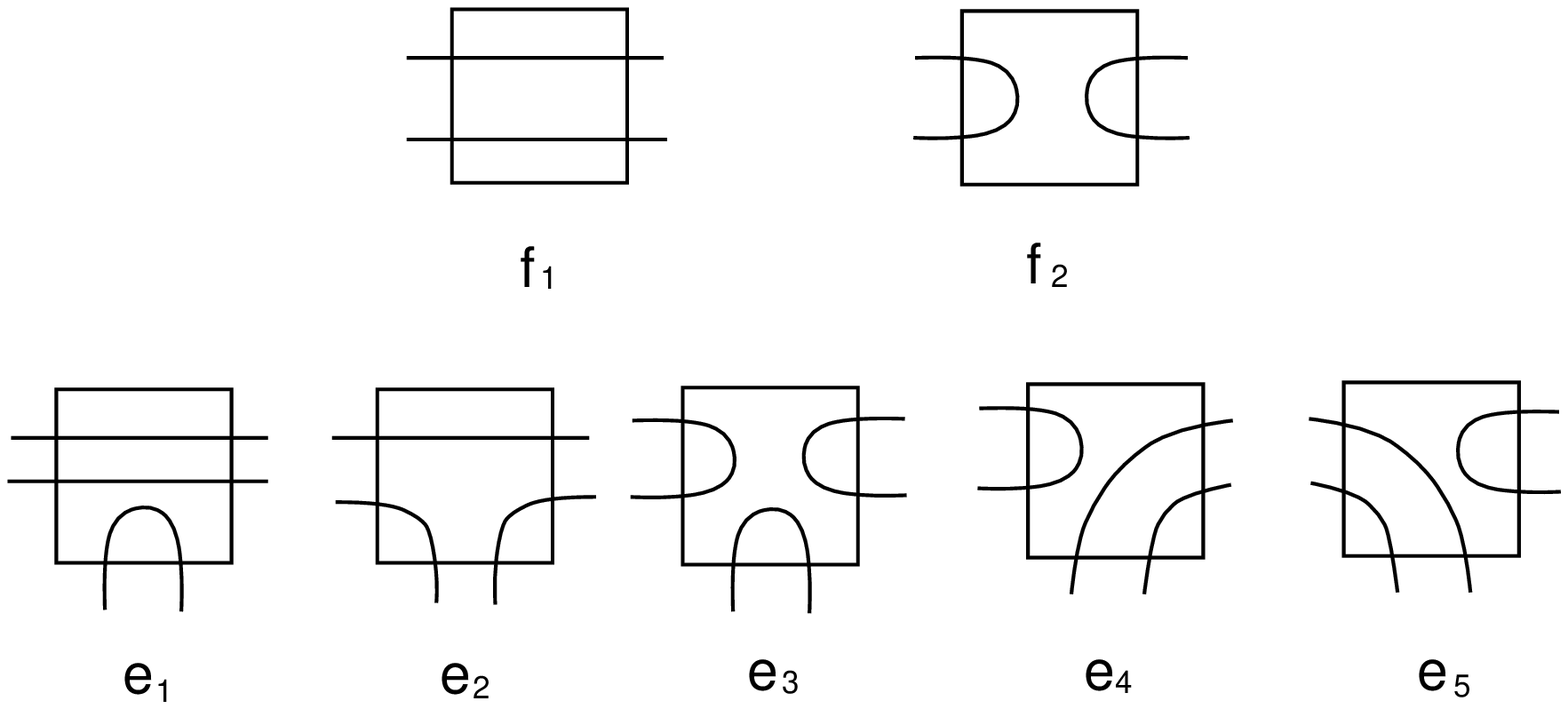}\\
\end{tabular}
\par\vspace{0.5cm}
Fig. 3.10
\end{center}

${\cal S}_{2,\infty}(3)$ is a left and a right ${\cal S}_{2,\infty}(2)$ module. 
The table for the right multiplication is shown in Figure 3.11, where
$\mu = -A^2-A^{-2}$.

\[
\left[ \begin{array}{cc}
e_1 & e_3 \\
e_2 & e_5 \\
e_3 & \mu e_3 \\
e_4 & e_3 \\
e_5 & \mu e_5
\end{array}
\right]
\]

\begin{center}     
Fig. 3.11. The $(i,j)$ entry of the matrix is equal to $(e_if_j)$.    
\end{center}   
  
Now let $L\in {\cal S}_{2,\infty}(3)$ and 
$L=a_1e_1 +a_2e_2 +a_3e_3 +a_4e_4 +a_5e_5$.
We look for $P=(xf_1 +yf_2) \in {\cal S}_{2,\infty}(2)$, such that 
$$ (**)\ \ \ \ \ \ \ \ \ \ \ \ \ \ L P=r_y (L P).$$
$L P = (a_1x)e_1 +(a_2x)e_2 +(a_1y+a_3x+\mu a_3y +a_4y)e_3 +(a_4x)e_4
+(a_2y +a_5x +\mu a_5y)e_5$. Notice that $r_y(f_i)= f_i$ and 
$r_y(e_i)=e_i$ for $i\leq 3$, and $r_y (e_4)=e_5,\ r_y(e_5)=e_4$.
Therefore (**) is equivalent to:
$a_4x=a_5x+a_2y+\mu a_5y$ or equivalently $y(a_2+\mu a_5)= x(a_4-a_5)$.
Now, either 
\begin{enumerate}
\item
[(i)] $a_2+\mu a_5=0$ and then $x(a_4-a_5)=0$, or
\item
[(ii)] 
 $a_2+\mu a_5\neq 0$ and then 
$\frac{y}{x} =\frac{a_4-a_5}{a_2+\mu a_5}$. Equivalently one has
 one projective solution $(x,y)=t(a_2+\mu a_5,a_4-a_5)$. 
\end{enumerate}

Now we have to check whether $P$ is invertible in ${\cal S}_{2,\infty}(2)$ 
(we are interested in two sided inverse so we use the fact that
${\cal S}_{2,\infty}(2)$ is commutative).
Let $Q=z_1f_1+z_2f_2$. Then $PQ=QP= (xz_1)f_1+ (yz_1+xz_2+\mu yz_2)f_2$.
Thus $Q$ is the inverse of $P$ iff $xz_1=1$ and $yz_1+ (x+\mu y)z_2=0$.
Thus $P$ is invertible iff $x\neq 0$ and $x+\mu y\neq 0$. 
If $P$ is invertible then 
$P^{-1}= \frac{1}{x}f_1+\frac{x-y}{(x+\mu y)x}f_2$.
From the above it follows that $L$ is in $D(3)$ iff $a_4=a_5$ or
$(a_2+\mu a_5)(a_2+\mu a_4)\neq 0$. Thus $D(3)$ contains the complement of
an algebraic set so it is dense in 
${\cal S}_{2,\infty}(3)$. 
\end{proof} 

\begin{theorem}\label{Theorem 3.4}
\begin{description}
\item
[(a)] Choose any tangles $L\in {\cal S}_{2,\infty}(3)$ and 
$T_i\in {\cal S}_{2,\infty}(2)$. Further choose any cyclic word, $w(L,T_i)$, 
over the alphabet $\{L,T_i\}$ and place the corresponding tangles
in the annulus as in Fig. 3.12(a). Now consider the cyclic 
word $w(r_z (L),T_i)$,
and again place the corresponding tangles in the annulus (Fig. 3.12(b)).
Then the elements of the Kauffman bracket skein module of the 
annulus corresponding to $w(L,T_i)$ and $w(r_z(L),T_i)$ are the same.

\par\vspace{1cm}
\begin{center}
\begin{tabular}{c}
\includegraphics[trim=0mm 0mm 0mm 0mm, width=.9\linewidth]
{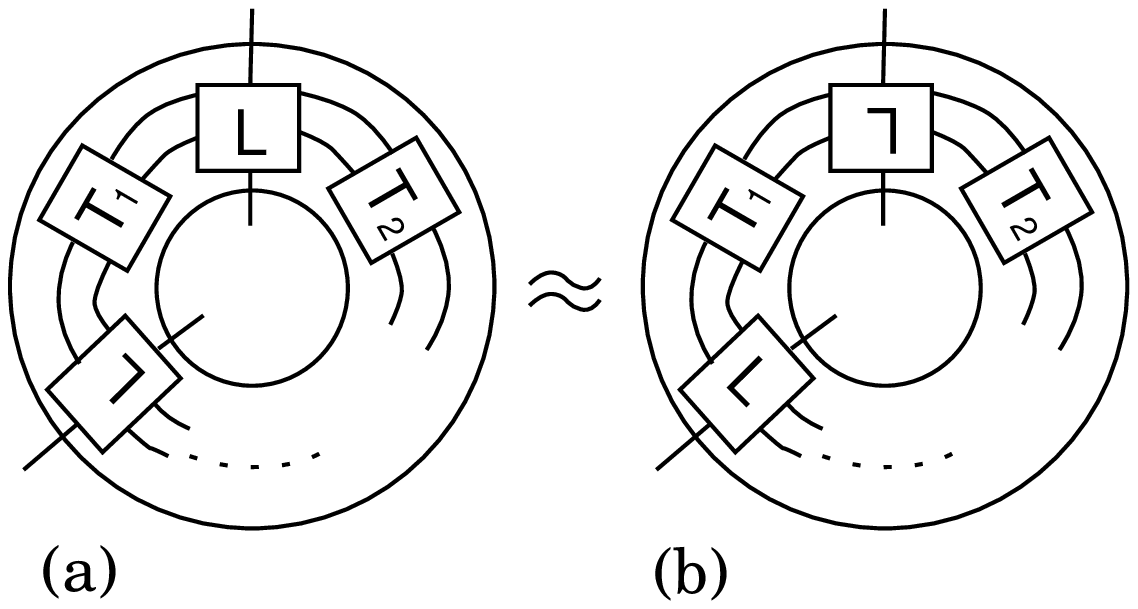}\\
\end{tabular}
\par\vspace{0.5cm}
Fig. 3.12
\end{center}

\item
[(b)]Choose any tangles $L\in {\cal S}_{2,\infty}(3)$ and tangles 
$T_i\in {\cal S}_{2,\infty}(2)$. Choose any cyclic word, $w(L,T_i)$, 
over the alphabet $\{L,T_i\}$ and place the corresponding tangles 
in the annulus as in Fig. 3.13(a). Now consider the cyclic 
word $w(r_y(L),T_i)$,
and again place the corresponding tangles in the annulus (Fig. 3.13(b)). 
Then the elements of the skein module of the annulus corresponding to
$w(L,T_i)$ and $w(r_y (L),T_i)$ are the same.

\par\vspace{1cm}
\begin{center}
\begin{tabular}{c}
\includegraphics[trim=0mm 0mm 0mm 0mm, width=.9\linewidth]
{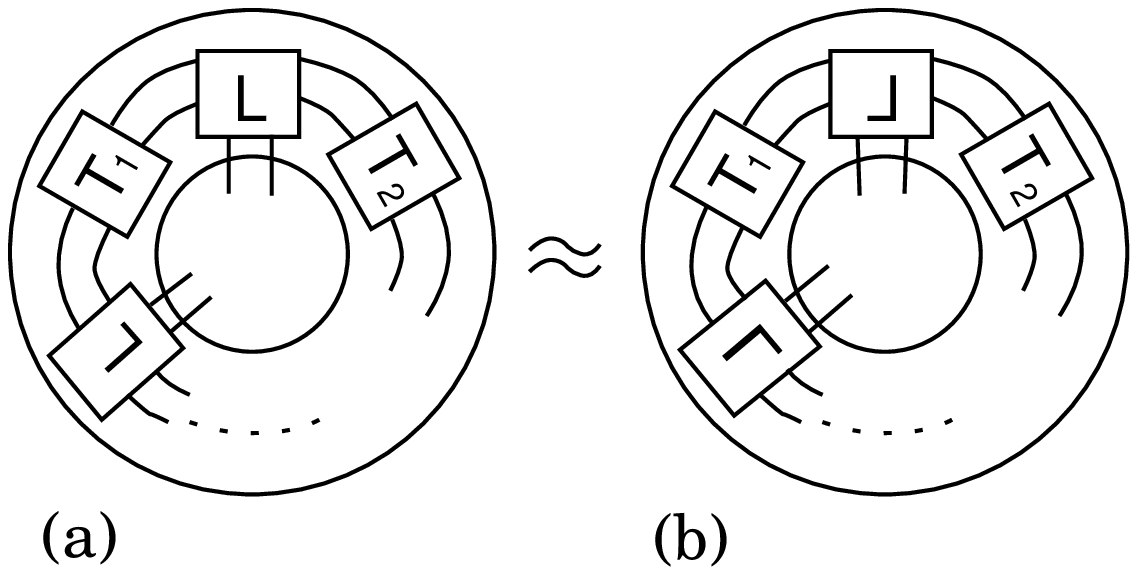}\\
\end{tabular}
\par\vspace{0.5cm}
Fig. 3.13
\end{center}

\end{description}
\end{theorem}
\begin{proof}
We will prove (b). The proof of (a) is analogous\footnote{
As in \cite{H-P-2}, we can deduce (a) from (b).}.
Let $L\in D(3)$ and $P\in {\cal S}_{2,\infty}(2)$ from Lemma 3.3(b).
Let us place $PP^{-1}$ in the annulus, as in Fig. 3.14. 
Then let $P$ travel
along the annulus, changing $L$ to $r_y (L)$, and finally canceling
$P$ with $P^{-1}$; thus Theorem 3.4 holds for a dense subset of
${\cal S}_{2,\infty}(3)$ and therefore for any element $L$ of
${\cal S}_{2,\infty}(3)$. 
This completes the proof of Theorem 3.4.
\end{proof}

\par\vspace{1cm}
\begin{center}
\begin{tabular}{c}
\includegraphics[trim=0mm 0mm 0mm 0mm, width=.3\linewidth]
{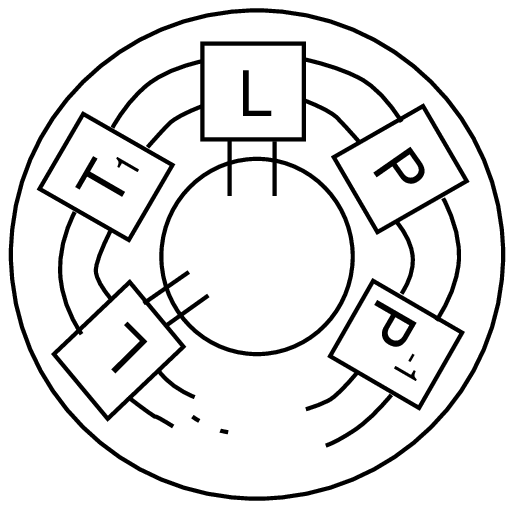}\\
\end{tabular}
\par\vspace{0.5cm}
Fig. 3.14
\end{center}

\section{Spectral parameter tangle for the skein polynomial}\label{4}
We describe in this section various applications of the Jones idea of
the spectral parameter tangle to oriented links, generalizing results of
Jones\cite{Jo-3} and Traczyk \cite{Tr-1}.
We can allow various 2-tangles and various mutations, for our 
construction. There are two essentially different ways of
orienting a 3-tangle: ``braid like" (Fig. 4.1(a)) and ``alternating"
(Fig. 4.1(b)). We concentrate here on the braid like case
because the Traczyk's method does not work in that case.

\par\vspace{1cm}
\begin{center}
\begin{tabular}{c}
\includegraphics[trim=0mm 0mm 0mm 0mm, width=.9\linewidth]
{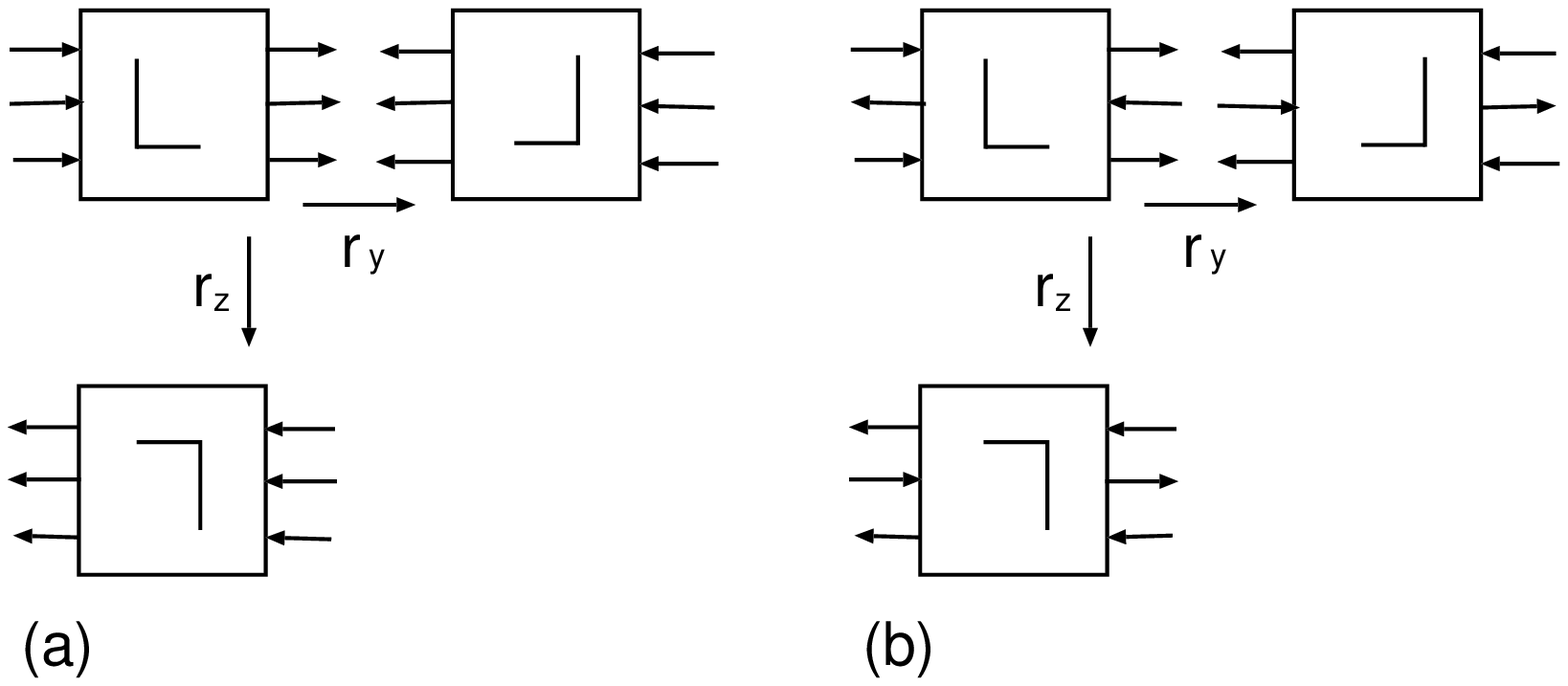}\\
\end{tabular}
\par\vspace{0.5cm}
Fig. 4.1
\end{center}

Recall that the skein (Homflypt) polynomial invariant of oriented 
links in $S^3$, $P_L \in Z[v^{\pm 1},z^{ \pm 1}]$, 
is given by:
\begin{description}
\item [(i)]
$P_{T_1} = 1$;
\item [(ii)]
$v^{-1}P_{L_+} - vP_{L_-} = zP_{L_0}.$
\end{description}
\par
The skein module (linear skein), ${\cal S}_3(M)$ is a generalization 
of the skein polynomial 
to any oriented 3-manifold (possibly with boundary). It is a quotient of
$R\cal L$ and the submodule generated by expressions
$v^{-1}P_{L_+} - vP_{L_-} - zP_{L_0}$, where $\cal L$ is the set of all
oriented links (including links with boundary on $\partial M$) 
up to the ambient isotopy, $R$ is commutative ring with unit containing
$Z[v^{\pm 1},z]$ and
$RX$ is the free $R$-module with basis $X$. Unless otherwise stated,
we work here with $R$ being the field of rational functions, ${\cal F}(v,z)$. 

Let ${\cal S}_3(n)$ denote the skein module of a tangle with
$n$ inputs and $n$ outputs as in Fig. 4.2 (inputs and outputs are fixed). 

\par\vspace{1cm}
\begin{center}
\begin{tabular}{c}
\includegraphics[trim=0mm 0mm 0mm 0mm, width=.3\linewidth]
{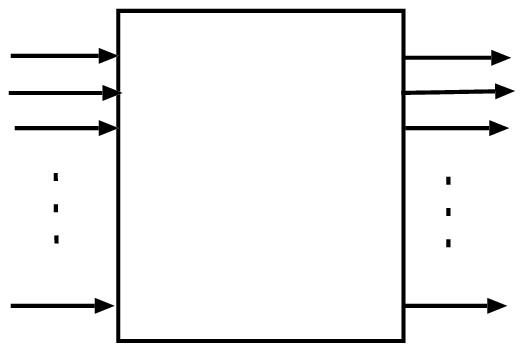}\\
\end{tabular}
\par\vspace{0.5cm}
Fig. 4.2
\end{center}

${\cal S}_3(n)$ is known to be the Hecke algebra of type $A$ and is free
with basis indexed by permutations, $S_n$ \cite{Bourbaki,Jo-1,M-T-2}.
In particular $(1,{\sigma}_1)$ is a basis of ${\cal S}_3(2)$ and
$(e_1,...,e_6)=(1,{\sigma}_1, {\sigma}_2, {\sigma}_1{\sigma}_2{\sigma}_1, 
{\sigma}_1{\sigma}_2, {\sigma}_2{\sigma}_1)$ is a basis of ${\cal S}_3(3)$; 
Fig. 4.3.

\par\vspace{1cm}
\begin{center}
\begin{tabular}{c}
\includegraphics[trim=0mm 0mm 0mm 0mm, width=.9\linewidth]
{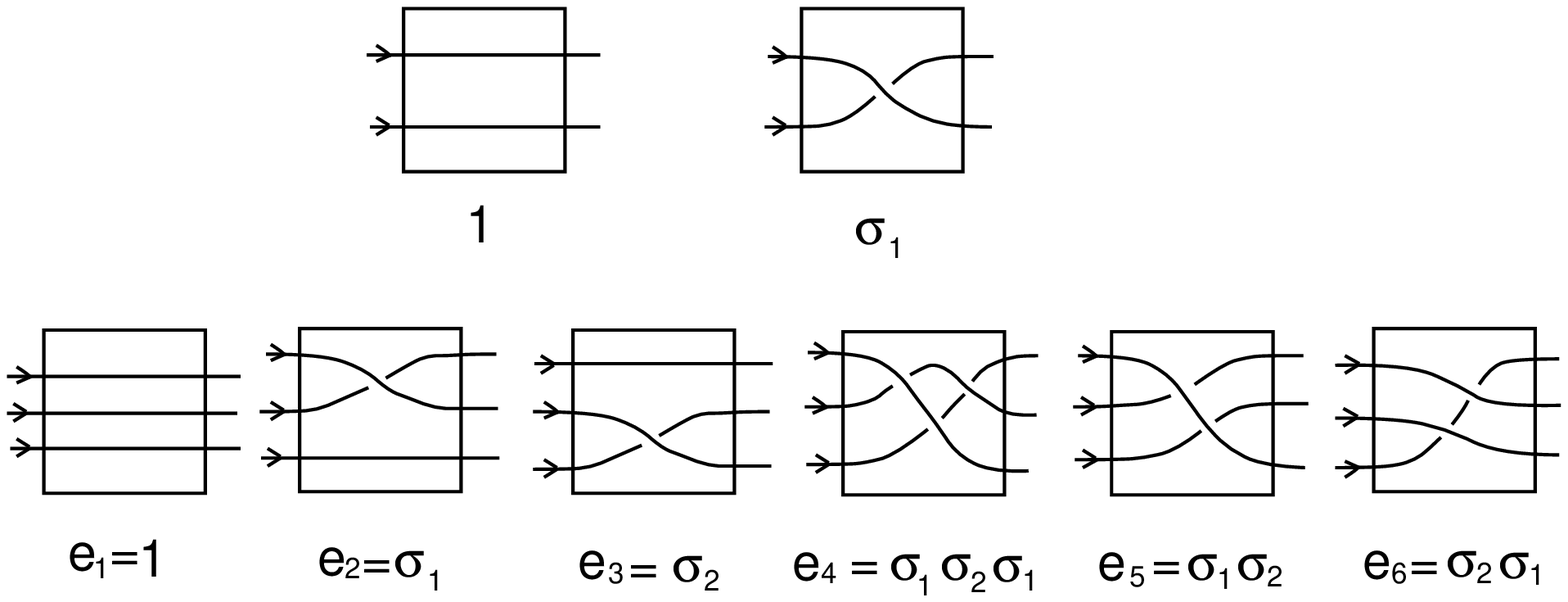}\\
\end{tabular}
\par\vspace{0.5cm}
Fig. 4.3
\end{center}

\begin{lemma}\label{Lemma 4.1}\
\begin{description}
\item
[(a)] 
Let $P= (x_1 + x_2{\sigma}_1)\in {\cal S}_3(2)$. 
$P$ is invertible in ${\cal S}_3(2)$ \  
iff \ \ \  
$x_1^2-v^2x_2^2 +vzx_1x_2\neq 0$.
\item
[(b)]
Let $L=(\Sigma_{i=1}^{6}a_ie_i)\in {\cal S}_3(3)$. If $LP={r_y}(LP)$ up
to the global change of orientation, then $x_1(a_5-a_6)=
x_2(a_3+vza_6-v^2a_4)$.
\item
[(c)]
Let $D_3(3)$ denote the subset of ${\cal S}_3(3)$ such that for 
$L \in D_3(3)$, there is an invertible element $P\in {\cal S}_3(2)$ 
satisfying  $LP=r_y (LP)$ (up
to the global change of orientation) in
${\cal S}_3(3)$ (see Fig. 4.4).  
Then $L\in D_3(3)$ iff $a_5=a_6$ or $F(a_3, a_4, a_5, a_6)\neq 0$ where
$F(a_3, a_4, a_5, a_6)=(a_3-v^2a_4+vza_5)(a_3-v^2a_4+vza_6)-v^2(a_5-a_6)^2$.

\par\vspace{1cm}
\begin{center}
\begin{tabular}{c}
\includegraphics[trim=0mm 0mm 0mm 0mm, width=.9\linewidth]
{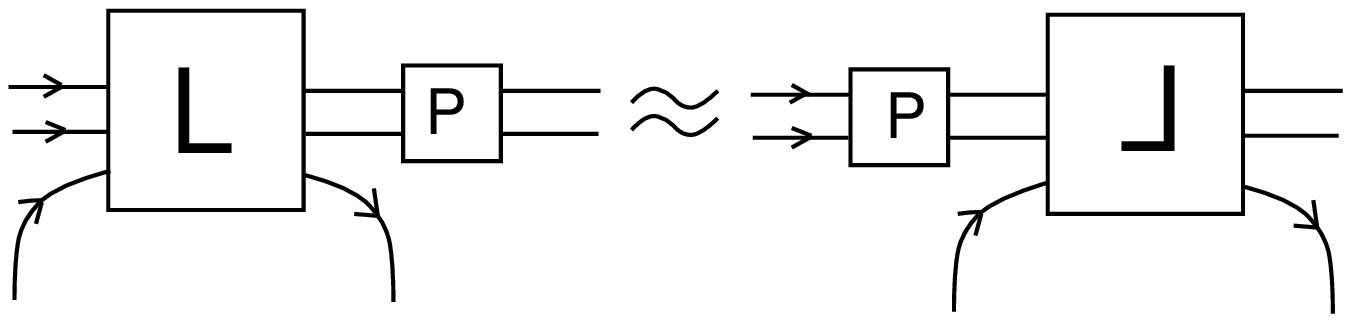}\\
\end{tabular}
\par\vspace{0.5cm}
Fig. 4.4
\end{center}

\item
[(d)] $D_3(3)$ is dense in ${\cal S}_3(3)$.
\end{description}
\end{lemma}
\begin{proof}
\begin{description}
\item
[(a)]
$P$ acts from the right side on ${\cal S}_3(2)$ and in the basis
$1,{\sigma}_1$ it is described by the matrix 
\[ [P] = \left[ \begin{array}{cc}
x_1 & v^2x_2 \\
x_2 & x_1+vzx_2
   \end{array}
\right]
\]
Thus $P$ is invertible iff 
det$[P] = x_1^2-v^2x_2^2 +vzx_1x_2\neq 0$.
\item
[(b)]
Let $<-,->$ be a bilinear form on ${\cal S}_3(3)$ given by
$<e_i,e_j>={\delta}_{i,j}$, then $LP={r_y}(LP)$ is equivalent to
$<LP,e_5>=<LP,e_6>$. The short calculation gives:
$<LP,e_5>=x_1a_5+x_2a_4v^2$, \ \ \ \ and $<LP,e_6>=x_1a_6+a_3x_2+vza_6x_2$. \
Thus $P$ is given by the condition $x_1(a_5-a_6)=x_2(a_3+vza_6-v^2a_4)$.

\item
[(c)]
If $a_5=a_6$ then we can take for example $P=P^{-1}=1$ otherwise we have
unique ``projective" solution 
$P=t((a_3+vza_6-v^2a_4) +(a_5-a_6)\sigma_1)$.
Putting this solution to the condition from (a) one gets:
$(a_3-v^2a_4+vza_5)(a_3-v^2a_4+vza_6)-v^2(a_5-a_6)^2\neq 0$.
\item
[(d)]
It is enough to notice that the last polynomial inequality holds for
a (open) dense subset of ${\cal S}_3(3)$ (a complement of an
algebraic set).
\end{description}
\end{proof}

\begin {theorem}\label{4.2}
\begin{enumerate}
\item
[(a)]  Let $L$ be any oriented 3-tangle such that 
the rotation $r_z$ sends inputs to outputs and vice versa, as shown
in Fig. 4.5(a). Let $o(L)$ denote the 3-tangle obtained from $L$ by changing
its orientation, Fig. 4.5(b). 

\par\vspace{1cm}
\begin{center}
\begin{tabular}{c}
\includegraphics[trim=0mm 0mm 0mm 0mm, width=.9\linewidth]
{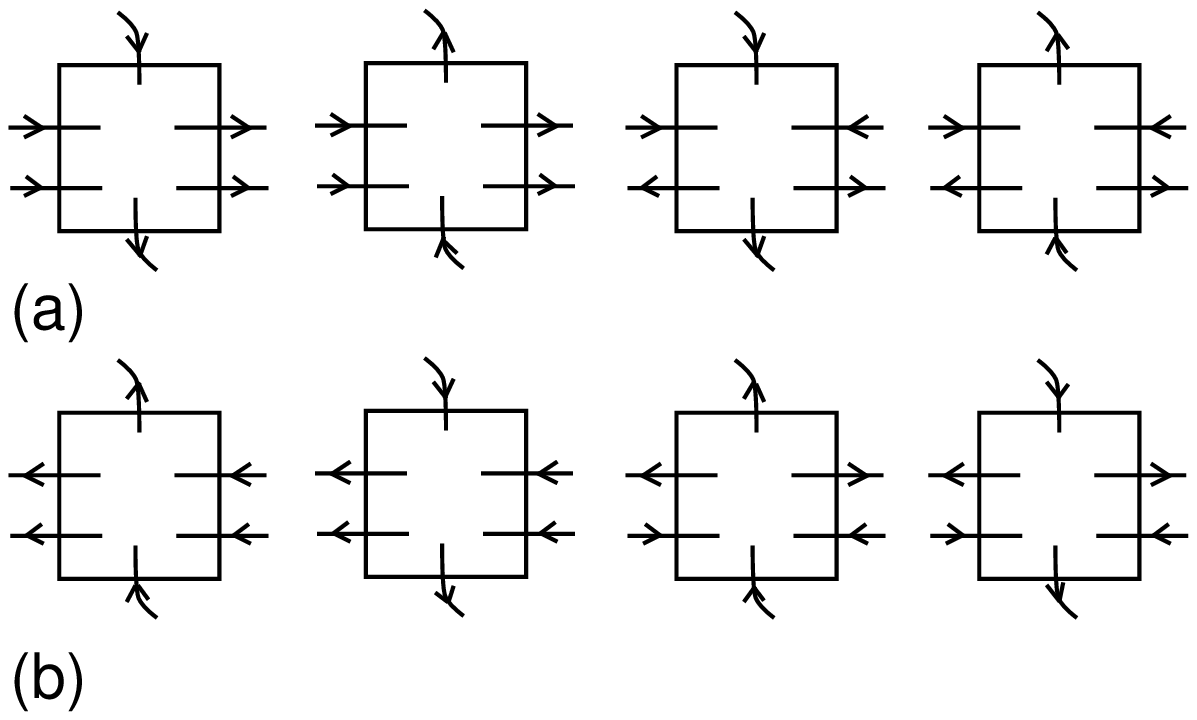}\\
\end{tabular}
\par\vspace{0.5cm}
Fig. 4.5
\end{center}

Consider now a cyclic 
word $w(L,o(L),D_{\alpha})$ over alphabet consisting of
tangles $L,o(L)$ and all oriented 2-tangles, $D_{\alpha}$. 
Furthermore assume that neighboring tangles in the cyclic word
have compatible orientation (can be glued together). Let $T(w)$ be 
the associated tangle placed in the annulus (oriented version
of Fig. 3.12(a)). Now, let us rotate,
along the axis $z$, each 3-tangle of the word, and then change its 
orientation.  We obtain, possibly different,
tangle but equal to the previous one in the skein module of the solid torus. 
That is
$T(w(L,o(L),D_{\alpha}))$ and 
$T(w(o(r_z(L)),o(r_z(o(L))),D_{\alpha}))$ are
equal in the skein module of the solid torus; compare Fig. 3.12. 
If the solid torus is embedded
in a 3-manifold and the endpoints of the tangles are connected in the 
same manner, outside the solid torus, then the resulting
links are equal in the skein module of the 3-manifold.
\item
[(b)]  Let $L$ be any oriented 3-tangle such that  
the rotation $r_y$ sends inputs to outputs and vice versa, as shown
in Fig. 4.6. 

\par\vspace{1cm}
\begin{center}
\begin{tabular}{c}
\includegraphics[trim=0mm 0mm 0mm 0mm, width=.9\linewidth]
{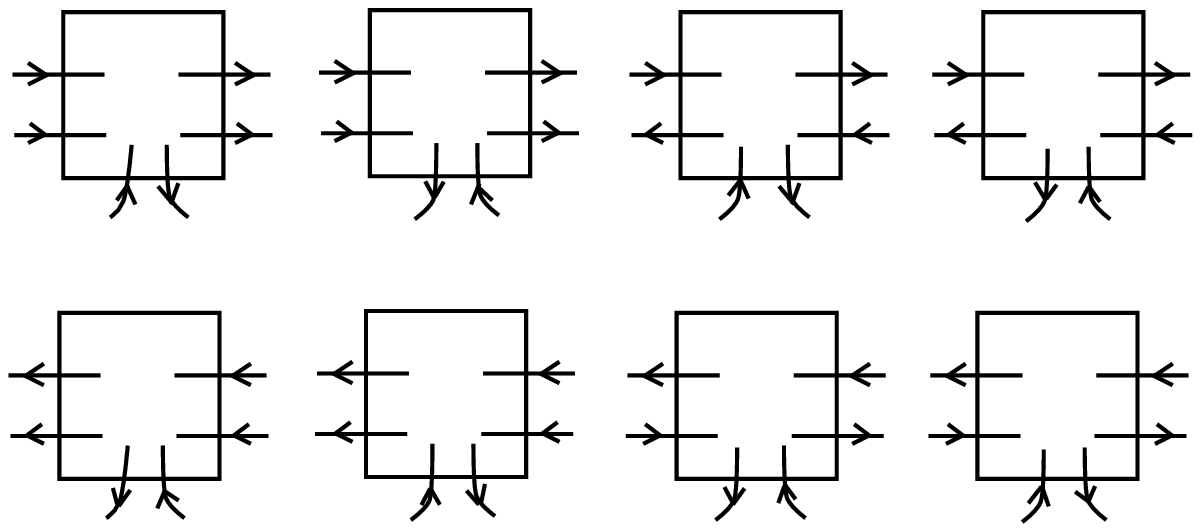}\\
\end{tabular}
\par\vspace{0.5cm}
Fig. 4.6
\end{center}

Consider now, a cyclic 
word $w(L,o(L),D_{\alpha})$ over the alphabet consisting of
tangles $L,o(L)$ and all 
oriented 2-tangles, $D_{\alpha}$. Furthermore assume that neighboring tangles in the cyclic word
have compatible orientation (can be glued together). Let $T(w)$ be the
associated tangle placed in the annulus. Now, let rotate,
along the axis $y$, each 3-tangle of the word,  and then change its 
orientation.  We obtain, possibly different,
tangle but equal to the previous one in the skein module of the solid 
torus. That is $T(w(L,o(L),D_{\alpha}))$ and\\
 $T(w(o(r_y(L)),o(r_y(o(L)),D_{\alpha})))$ are
equal in the skein module of the solid torus; compare Fig. 3.13. 
If the solid torus is embedded
in a 3-manifold and endpoints of the tangles are connected in the 
same manner, outside the solid torus, then the resulting
links are equal in the skein module of the 3-manifold.
\end{enumerate}
\end{theorem}

\begin{proof}
If $L$ is a ``braid like" 3-tangle (Fig. 4.1(a)), we use the spectral 
parameter tangle (from Lemma 4.1) similarly as in the
case of the Kauffman bracket skein relation. All cases of the orientation
of $L$ from Fig. 4.5 and 4.6 easily reduce to the basic cases of
``braid like" and ``alternating" 3-tangles. For an ``alternating"
3-tangles one should prove a lemma analogous to Lemma 4.1, but there is
no difficulty in doing so.
\end{proof}

\section{Spectral parameter 3-tangle }\label{5}
We use in this section the idea of the spectral parameter tangle but 
in the more involved case. The spectral parameter tangle is, 
in this section, a 3-tangle. 
We work in this part with the third (Homflypt) skein relation, i.e. with
the skein module ${\cal S}_3$. For the Kauffman bracket skein module,
computations are similar but slightly shorter. I delayed publishing
these results, which were ready in the summer of 1991 \cite{P-ab}, 
because I believed that a
similar result could hold for the Kauffman polynomial. Only two years
later Traczyk performed calculations which showed that it is not a case
(even for a 2-cable of the Jones polynomial which gives a special substitution
of the Kauffman polynomial). 
\begin {theorem}\label{5.1}\ \\ 
If $X$ and $Y$ are 3-tangles oriented as in Figure 5.1(a) (resp. Fig. 5.1(b))
 and $W(XY)$ is any word in letters $X$ and $Y$ then
$$Tr(W(X,Y))=Tr(W(r_y(X),r_y(Y))) \  in\   {\cal S}_3(S^1 \times D^2)$$ \\
where
 $ Tr(A)$ denotes the cyclic closure of the $3$-tangle $A$  (i.e. the link 
diagram in the annulus determined by $A$); compare Fig. 5.2.
\end{theorem}

\par\vspace{1cm}
\begin{center}
\begin{tabular}{c}
\includegraphics[trim=0mm 0mm 0mm 0mm, width=.5\linewidth]
{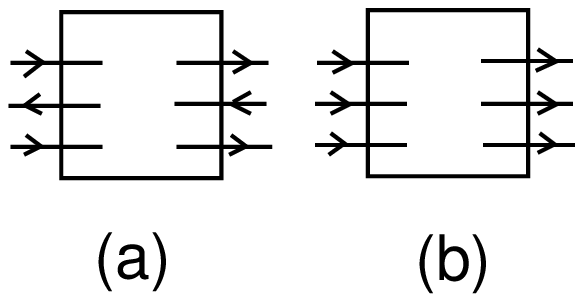}\\
\end{tabular}
\par\vspace{0.5cm}
Fig. 5.1
\end{center}

Theorem 5.1 contains 2 parts depending on  whether tangles
$X$ and $Y$ are of type (a) or (b) of Figure 5.1. 
The idea in both cases is the same
but specific calculations different. We  prove here the simpler case
of tangles $X,Y$ oriented as in Figure 5.1(a).\footnote{ 
The initial proof of the second case with segments oriented as 
in Figure 5.1(b) was performed with
the using the computer program {\it Mathematica} (and with help of J. Walsh). 
We will present a ``computer free" proof in the future paper \cite{P-4}.}
As in previous parts the crucial element of the proof is the existence
of the invertible spectral parameter tangle and then Theorem 5.1 follows
easily, as the results in the previous parts. Denote by ${\cal S}'_3(3)$ the
skein module of the tangle from Fig. 5.1(a).

\par\vspace{1cm}
\begin{center}
\begin{tabular}{c}
\includegraphics[trim=0mm 0mm 0mm 0mm, width=.9\linewidth]
{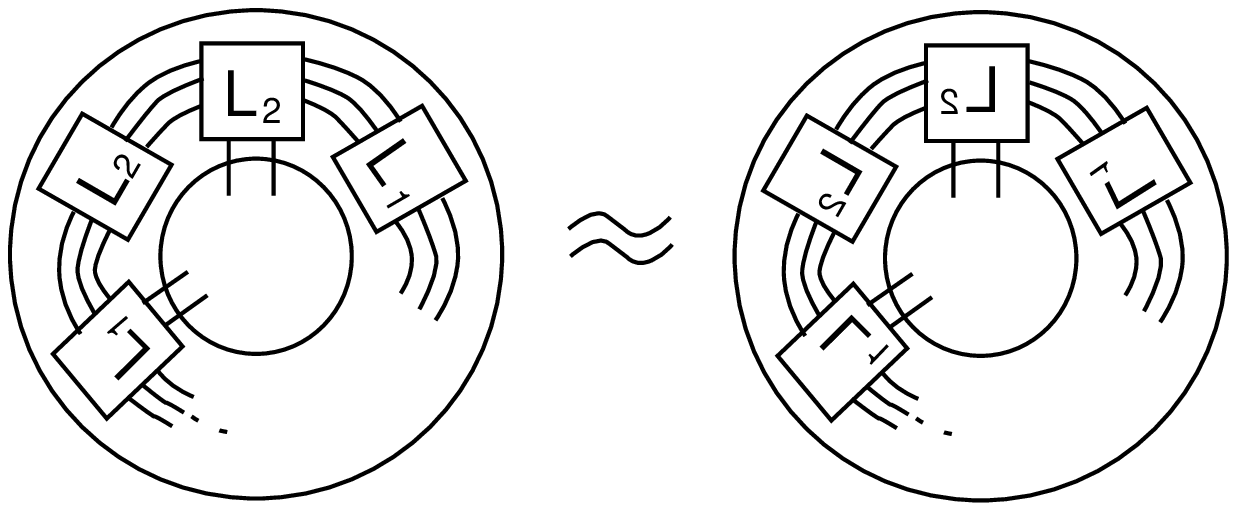}\\
\end{tabular}
\par\vspace{0.5cm}
Fig. 5.2
\end{center}

\begin{lemma}\label{5.2}
\begin {enumerate}
\item
[(a)] Let $Y={\Sigma}_{i=0}^{5}y_ie_i \in {\cal S}'_3(3)$, in the 
basis of ${\cal S'}_3(3)$ shown in Fig. 5.3.\ $Y$ is invertible
in ${\cal S}'_3(3)$ \ iff \ 
\begin{eqnarray*}
((y_0+y_3+\mu y_1)(y_0+y_4+\mu y_2)&-&(y_2+y_5+\mu y_4)(y_1+y_5+\mu y_3))\\
&\times&(y_0(y_0+v^{-1}zy_5)-v^{-2}y_5^2)\neq 0.
\end{eqnarray*}
\item
[(b)] Let $X={\Sigma}_{i=0}^{5}x_ie_i \in {\cal S}'_3(3)$.
If $XY=r_y(XY)$ and $Y=r_y (Y)$ (up to the global change of orientation)
 then $y_3=y_4$ and
$y_0(x_4-x_3)+y_1(x_2+x_5+\mu x_4)-y_2(x_1+x_5\mu x_3)$\\
$-y_3(x_3-x_4 +\mu (x_1-x_2))+y_5(x_2-x_1) = 0 .$
\item
[(c)]
Let $D^2(3)$ denote the subset of ${\cal S}'_3(3)\times {\cal S}'_3(3)$,
such that $(A,B)\in D^2(3)$ iff there is an invertible element $Y\in
{\cal S}'_3(3)$ such that $AY=Yr_y(A)$ and $BY=Yr_y(B)$. Then
$D^2(3)$ is dense in ${\cal S}'_3(3)\times {\cal S}'_3(3)$.
\end{enumerate}
\end{lemma}
~~
\par\noindent
\begin{proof}
\begin{enumerate}
\item
[(a)]
Let $X={\Sigma}_{i=0}^{5}x_ie_i$,  $Y={\Sigma}_{i=0}^{5}y_ie_i$, in
the basis of ${\cal S}'_3(3)$ shown in Fig. 5.3. 

\par\vspace{1cm}
\begin{center}
\begin{tabular}{c}
\includegraphics[trim=0mm 0mm 0mm 0mm, width=1.0\linewidth]
{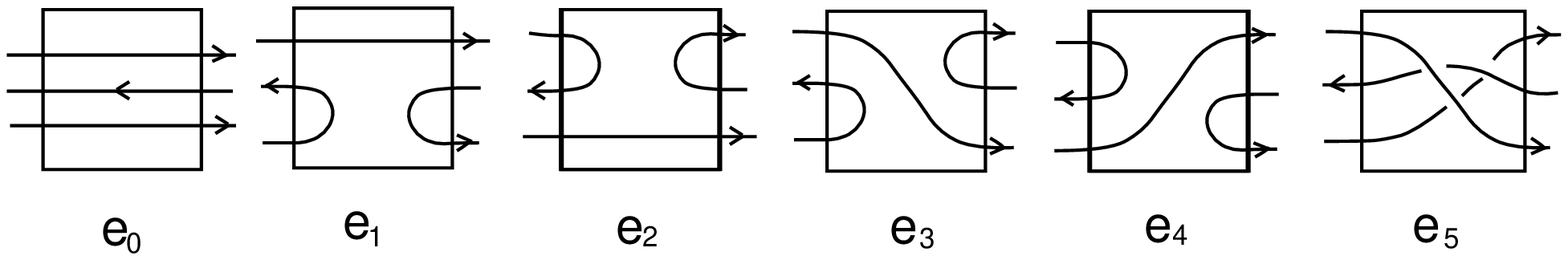}\\
\end{tabular}
\par\vspace{0.5cm}
Fig. 5.3
\end{center}
\par\vspace{1cm}

Then $XY=$
$$e_0 (y_0x_0+v^{-2}y_5x_5)+$$
$$e_1 (y_0x_1+y_1(x_0+x_3+\mu x_1)+\ \ y_4(x_1+x_5+\mu x_3)+
y_5(x_3-v^{-1}zx_5))+$$
$$e_2(y_0x_2+\ \ y_2 (x_0+x_4+\mu x_2)+y_3(x_2+x_5+\mu x_4) +
y_5(x_4 -v^{-1}zx_5))+$$
$$e_3(y_0x_3+\ \ y_2(x_1+x_5+\mu x_3)+y_3(x_0+x_3+\mu x_1)+\ \ 
y_5x_1)+$$
$$e_4(y_0x_4+y_1(x_2+x_5+\mu x_4) +\ \ \  \ y_4(x_0+x_4+\mu x_2)+
y_5x_2)+$$
$$e_5(y_0x_5+ \ \ \ \ \ \ \ \ \ \ \ \ \ y_5(x_0+v^{-1}zx_5))$$
If $X$ is treated as a linear operator acting (from the left) on 
${\cal S}'_3(3)$,
then from the above calculation we obtain that in the basis $e_0,e_1,...,e_5$:
\begin{eqnarray*}
\det X&=&((x_0+x_3+\mu x_1)(x_0+x_4+\mu x_2)\cr
&-&(x_2+x_5+\mu x_4)(x_1+x_5+\mu x_3))^2(x_0(x_0+v^{-1}zx_5)-v^{-2}x_5^2).
\end{eqnarray*}
Then the part (a) follows (we have to change roles of $X$ and $Y$).
\item
[(b)] follows immediately from the formula for $XY$.
\item 
[(c)] One have to perform easy but long and tedious calculations.
\end{enumerate}
\end{proof}

One of the ingredients of our work was to show that the elements of the 
Temperley-Lieb algebra $TL_2$ and $TL_3$ are almost 
everywhere invertible. We can show that it holds for any $n$, 
answering the question raised by Rolfsen in \cite{Ro} (Question after
Proposition 5). Similar fact 
holds also for the Hecke algebra $H_n(v,z)={\cal S}_3(n)$ as well
as for Birman-Murakami-Wenzl algebra.
\par~~\par
\begin{theorem}\label{5.3}\ \\
\begin{enumerate}
\item
[(a)] Invertible elements of $TL_n={\cal S}_{2,\infty}(n)$ form
a dense subset of $TL_n$.
\item
[(b)] Invertible elements of $H_n(v,z)$ form
a dense subset of $H_n(v,z)$. 
\end{enumerate}
\end{theorem}
\begin{proof}
We will demonstrate (b). The proof of (a) is similar. The main idea of
the proof is observation that the analogous result for the group
algebra over the symmetric group $S_n$ can be seen immediately\footnote
{We could stop here saying that $H_n(v,z)$ and ${\cal F}(v,z)S_n$ are
algebra isomorphic, but we present an elementary proof.}.

It is convenient to work with the ring $Q[v^{\pm 1},z]$, in 
addition to the field ${\cal F}(v,z)$. Denote the corresponding 
(to the ring $Q[v^{\pm 1},z]$) skein
module and Hecke algebra by ${\cal S}^+_3(n)=H^+_n(v,z)$.
Let $\{e_1=1,e_2,...,e_{n!}\}$ be the basis of $H_n(v,z)$, composed 
of positive
braids of minimal number of crossings for a given permutation.
It is also a basis of the free module $H^+_n(v,z)$.
Let $X$ be an element of $H_n(v,z)$ and consider its action from the left
on $H_n(v,z)$. It is a linear function and its matrix $[X]=[a_{i,j}]$ in our
basis is given by $Xe_i={\Sigma}_ja_{i,j}e_j$.
Invertible elements of $H_n(v,z)$ form a (multiplicative) 
subgroup of $H_n(v,z)$. It is
immediate but important observation that in order to show that they are
dense in $H_n(v,z)$, it suffices to show that elements of non-zero
determinant in $H^+_n(v,z)$ are dense in $H^+_n(v,z)$. Consider
$X={\Sigma}_{i=1}^{n!}x_ie_i$ in $H^+_n(v,z)$. Because $x_i \in Q[v^{\pm 1},z]$,
we can work now modulo $(z,v^2-1)$ -- the ideal generated by $z$ and $v^2-1$.
Our skein relation ${\sigma}_i^2=vz{\sigma}_i+v^2$ reduces to
${\sigma}_i^2\equiv 1\ mod(z,v^2-1)$. We can immediately check that
$det X \equiv x_1^{n!}+ {\cal O}(x_1^{n!-2}) mod(z,v^2-1)$, because for
$[X]=[a_{i,j}]$, $a_{i,i}\equiv x_1 mod(z,v^2-1)$ and $a_{i,j}\equiv
x_s mod(z,v^2-1)$ where $s \neq 1$ for $i\neq j$. Thus $det X \neq 0$
outside an algebraic set and Theorem 5.3(b) holds.
\end{proof}

Another application of the spectral parameter tangle has been 
recently demonstrated by Rolfsen in \cite{Ro-1}. Furthermore, 
the remarkable example by T.Kanenobu of an infinite family 
of different knots with the same skein polynomial \cite{Kan-1} 
can be put in more general
context similar to that of the spectral parameter tangle \cite{Kan-2}.
Namely Kanenobu considers elements of ${\cal S}'_3(n)$ for which $e_2$ of
Figure 5.3 is an eigenvector.

\section{Streszczenie}
Dok\l adnie 10 lat temu\footnote{Added for e-print: Summary (Streszczenie 
in Polish) is taken from my Banach Center presentation, March 24, 1994.}, 
wiosn\c a 1984, Vaughan Jones skonstruowa{\l} 
nowy, bardzo skuteczny niezmiennik w\c ez\l\'ow. Szybko okaza\l o si\c e,
\.ze wielomian Jonesa nie zawsze odr\'o\.znia r\'o\.zne w\c ez\l y.
Pozosta\l o jednak otwartym pytanie czy nietrywialny w\c eze\l \ mo\.ze
mie\'c trywialny wielomian Jonesa. Om\'owimy, na tym wyk\l adzie, trzy
podej\'scia do tego problemu, wszystkie uog\'olniaj{\c a}ce 
konstrukcj{\c e} {\it mutacji} Conwaya.
\begin{enumerate}
\item [1.] 
{\L}{\c a}cz{\c a}c ide{\c e} mutacji i satelity splotu.
\item [2.]
U\.zywaj\c ac idei rotoru. Pomys\l \ ten zaczerpni\c ety jest z teorii
graf\'ow.
\item
[3.]
 U\.zywaj\c ac idei sup\l a ze spektralnym parametrem. Jest to pomys{\l} 
Jonesa, bazuj\c acy na idei  Baxtera (stosowanej w teorii dok\l adnie
rozwi\c azywalnych modeli mechaniki statystycznej).
\end{enumerate}

J\'ozef H.~Przytycki\\
Department of Mathematics\\
and Computer Science\\
Odense University\\
Denmark\\
e-mail: $Jozef@imada.ou.dk$

and
\ \\
Warsaw University, Poland\\
e-mail: $Jozef@mimuw.edu.pl$\\
\ \\
Current address (Sep. 94 - Jun. 95):\\
Department of Mathematics\\
University of California\\
Berkeley, CA 94 720\\
USA\\
e-mail: $Jozef@math.berkeley.edu$\\
\ \\
(Added for e-print:\ \  from August 1995:\\
Department of Mathematics\\
George Washington University,\\
Washington, DC\\
e-mail: $przytyck@gwu.edu$.)


\begin{thebibliography}{99}

\bibitem
{Alb}
D.Albers, John Horton Conway -- Talking a good game, {\em Math. Horizons}, 
Spring 1994, Published by the M.A.A.

\bibitem
{A-P-R}
R.P.Anstee, J.H.Przytycki, D.Rolfsen. Knot polynomials and
generalized mutation.  {\em Topology and Applications}
32 (1989) 237-249.\ (Added for e-print:\\
http://front.math.ucdavis.edu/math.GT/0405382 )

 
\bibitem
{Baxter}
R.Baxter, {\em Exactly solved models in statistical mechanics},
Academic Press, London, 1982.

\bibitem
{Bourbaki}
N.Bourbaki, {\em Groupes et alg\`ebres de Lie, VI: Groupes de 
Coxeter et syst\`emes de Tits}, Herman paris, 1968.

\bibitem
{BSST} R.I.Brooks, C.A.B.Smith, A.H.Stone, W.T.Tutte. The Dissection of 
rectangles into squares. {\em Duke Math. Jour.}, 7 (1940) 312 -340.

\bibitem
{Co}
J.~H.~Conway, An enumeration of knots and links, {\em Computational
problems in abstract algebra} (ed. J.Leech),
Pergamon Press (1969) 329 - 358.

\bibitem
{HOMFLY}
P.~Freyd, D.~Yetter, J.~Hoste, W.~B.~R.~Lickorish,
K.~Millett, A.~Ocneanu, A new polynomial invariant of knots
and links, {\em Bull. Amer. Math. Soc.}, 12 (1985) 239-249.

\bibitem
{Hos} J.Hoste, A polynomial invariant of knots and links.
{\em Pacific J. Math.}, 124 (1986) 295-320.

\bibitem
{H-P-1} J.Hoste, J.H.~Przytycki, A survey of skein modules of 
3-manifolds,
``Knots 90", De Gruyter, Berlin, New York, 1992, 363-379.

\bibitem
{H-P-2}
J.Hoste, J.H.~Przytycki,  Tangle surgeries which preserve Jones-type
polynomials, Center for Pure and Applied Mathematics preprint-PAM 617, 
U.C.Berkeley, 1994. (Added for e-print: 
{\it International Journal of Mathematics} 8, 1997, 1015--1027.)

\bibitem
{J-R}
G.T.Jin, D.Rolfsen, Some remarks on rotors in link theory,
{\em Canadian
Math. Bull.}, 34 (1991), 480-484.

\bibitem
{Jo-1} V.~Jones. Hecke algebra representations of braid groups
and link polynomials, {\em Ann. of Math.}, 126 (1987) 335-388.

\bibitem
{Jo-2}
V.F.R.~Jones, On knot invariants related to some statistical mechanical
models, {\em Pacific J. Math.}, 137 (1989), 311-334.

\bibitem
{Jo-conf}
V.F.R.~Jones, Talk given at 25th Annual Spring Topology Conference, CSU
Sacramento, April 11, 1991.

\bibitem
{Jo-3} V.F.R.~Jones, Commuting transfer matrices and link polynomials,
{\em International Journal of Math.} 3(1992), 205-212.

\bibitem
{Jo-4} V.F.R.~Jones, Coincident link polynomials from
commuting transfer matrices, {\em Proceedings of the
XXth International Conference on Differential Geometric
Methods in Theoretical Physics}, Vol.1,2 (New York, 1991),
137-151, World Sci. Publishing, River Edge, NJ, 1992.

\bibitem
{Kan-1} T.Kanenobu, Infinitely many knots with the same polynomial invariant,
{\em Proc. Amer. Math. Soc.}, 97(1), (1986), 158-162.

\bibitem
{Kan-2} T.Kanenobu, The Homfly and the Kauffman bracket polynomials 
for the generalized mutant of a link, {\em Topology and its Applications},
to appear. (Added for e-print: 61(3), 1995, 257--279.)

\bibitem
{Kania}
J. Kania-Bartoszy\'nska, Examples of different 3-manifolds with 
the same invariants of Witten and Reshetikhin-Turaev, {\em Topology}, 
32, (1993), 47-54.

\bibitem
{Ka-1}
L.H.~Kauffman, State models and the Jones polynomial,
{\em Topology 26} (1987) 395-407.

\bibitem {L-bull} W.B.R.Lickorish, Polynomials for links, 
{\em  Bull. London Math. Soc.}, 20 (1988) 558-588.

\bibitem {Li-1}
W.B.R.Lickorish, Distinct 3-manifolds with all $SU(2)_q$ invariants
the same, {\it Proc. Amer. Math. Soc.}, 117 (1993), 285-292. 

\bibitem
{L-L} 
W.B.R.~Lickorish, A.S.~Lipson, Polynomials of 2-cable-like links,
{\em Proc. Amer. Math. Soc.}, 100 (1987), 355-361.

\bibitem
{L-M}  W.B.R.~Lickorish, K.~Millett, A polynomial invariant
of oriented links, {\em Topology} 26 (1987), 107-141.



\bibitem
{Mon} J.M.Montesinos, Surgery on links and double branched covers of $S^3$;
in {\em Knots, groups and 3-manifolds}, ed. L.P.Neuwirth, Ann. Math. 
Studies, 84, 227-259, Princeton Univ. Press, 1975.
 
\bibitem {M-T-1} H.R.Morton, P.Traczyk,  The Jones polynomial of
satellite links around mutants,
in {\em Braids}, Ed. J.S.Birman, A.Libgober,
AMS Contemporary Math., 78(1988), 587-592.

\bibitem
{M-T-2}
H.R.Morton, P.Traczyk, Knots and algebras, {\em Contribuciones Matematicas
en homenaje al profesor D.Antonio Plans Sanz de Bremond},
ed. E.Martin-Peinador and A.Rodez Usan, University of Zaragoza,
(1990), 201-220.

\bibitem {Mura}
J.Murakami, The parallel version of polynomial invariants of links,
{\em Osaka J. Math.}, 26(1989), 1-55.

\bibitem
{P-1}
J.H.~Przytycki, Equivalence of cables of mutants of knots,
{\em Canad. J. Math.}, XLI(2), (1989), 250-273.

\bibitem
{P-2} 
J.H.~Przytycki, Skein modules of 3-manifolds, {\em Bull. Ac. Pol.: Math.}, 
39(1-2), (1991), 91-100.

\bibitem
{P-3}
J.H.~Przytycki, Manuscript of the lecture delivered at the University
of Tennessee, October 18, (1991).

\bibitem
{P-ab}
J.H.~Przytycki, Applications of the spectral parameter tangle of V.~Jones,
{\em Abstracts of AMS}, 12 (1991), 496-497.

\bibitem
{P-4}
J.H.~Przytycki, The spectral parameter 3-string tangle, in preparation.

\bibitem
{PT} J.~H.~Przytycki, P.~Traczyk, Invariants of links of Conway type,
{\em Kobe J. Math. } 4 (1987) 115-139.

\bibitem
{R-T} N.Y.Reshetikhin, V.Turaev, Invariants of three manifolds via link
polynomials and quantum groups. {\em Invent. Math. } 103 (1991) 547-597.

\bibitem
{Ro} D.~Rolfsen, The quest for a knot with trivial Jones polynomial;
diagram surgery and the Temperley-Lieb algebra, in :{\em Topics in
knot theory}, Ed. M.E.Bozh{\"u}y{\"u}k, NATO ASI Series, Series C:
Mathematical and Physical Sciences - Vol. 399, Kluwer Academic
Publishers (1993), 195-210. 

\bibitem
{Ro-1} D.~Rolfsen, Global mutation of knots, {\em Journal of Knot
Theory and its Ramifications}, 3(3) (1994), 407-417.

\bibitem
{T-L} H.N.V.Temperley, E.H. Lieb. Relations between the ``percolation" 
and ``coloring" problem and other graph-theoretical problems
associated with regular planar lattices:
some exact results for the ``percolation" problem,
{\em Proc. Roy. Soc. }, London Ser. A 322 (1971), 251-280.

\bibitem
{Tr-1}
P.Traczyk, A note on rotant links, preprint 1989.
(Added for e-print: {\it J. Knot Theory Ramifications}, 3 (1994), 407-417.)

\bibitem
{Tu}
V.Turaev, The Yang-Baxter equation and invariants of links, {\em Invent.
Math.}, 92 (1988), 527-553.

\bibitem
{Tu-1} V.G.~Turaev, The Conway and Kauffman modules of the
solid torus, {\em Zap. Nauchn. Sem. Lomi} 167 (1988), 79-89. 
(Added for e-print: English translation: {\it J. Soviet Math.}, 
52, 1990, 2799-2805).

\bibitem
{T-W}
V.G.~Turaev, H.Wenzl, Quantum invariants of 3-manifolds associated with
classical simple Lie algebras, {\em Internat. J. Math.}, 4(2), (1993), 323-358.

\bibitem{Tut}
W.T.Tutte, Codichromatic graphs, {\em J. Combin. Theory}, Ser. B, 16 (1974),
168-174.

\bibitem
{Viro} O.Ya.Viro, Nonprojecting isotopies and knots with homeomorphic 
coverings, {\it Jour. Sov. Math.}, 12, (1979), 86-96. 

\bibitem
{Wa}
F.Waldhausen, {\"U}ber Involutionen der 3 Sph\"are,
{\em Topology}, 8, (1969), 81-91.

\bibitem
{Ya}
S.Yamada, An operator on regular isotopy invariants of
link diagrams, {\em Topology}, 28(3), (1989), 369-377.

\end{thebibliography}
\end{document}